\documentclass[UTF-8,reqno]{amsart}
\usepackage{enumerate}
\usepackage{amssymb,url,color, booktabs}
 \usepackage{mathrsfs}
 
\usepackage{amssymb}
\usepackage{txfonts}
\usepackage{bbm}
\usepackage{cases}
\usepackage{amsmath, bbm}
\usepackage{graphicx}
\usepackage{mathrsfs}
\usepackage{stmaryrd}
\usepackage{amsfonts}
\usepackage{enumerate,amsmath,amssymb,amsthm}

\usepackage[colorlinks,linkcolor=black,citecolor=blue]{hyperref}

\numberwithin{equation}{section}

\newcommand{\be}{\begin{eqnarray}}
\newcommand{\ee}{\end{eqnarray}}
\newcommand{\ce}{\begin{eqnarray*}}
\newcommand{\de}{\end{eqnarray*}}
\newtheorem{theorem}{Theorem}[section]
\newtheorem{lemma}[theorem]{Lemma}
\newtheorem{remark}[theorem]{Remark}
\newtheorem{definition}[theorem]{Definition}
\newtheorem{proposition}[theorem]{Proposition}
\newtheorem{Examples}[theorem]{Example}
\newtheorem{corollary}[theorem]{Corollary}

\def\1{{\bf 1}}

\def\eps{\varepsilon}

\def\e{\mathrm{e}}

\def\p{\partial}

\def\[{{\Big[}}
\def\]{{\Big]}}
\def\<{{\langle}}
\def\>{{\rangle}}
\def\({{\Big(}}
\def\){{\Big)}}

\def\bx{{\mathbf{x}}}

\def\dif{{\mathord{{\rm d}}}}

\def\min{{\mathord{{\rm min}}}}
\def\Vol{\mathord{{\rm Vol}}}

\def\no{\nonumber}
\def\={&\!\!=\!\!&}

\def\cA{{\mathcal A}}
\def\cB{{\mathcal B}}

\def\cJ{{\mathcal J}}

\def\mD{{\mathbb D}}
\def\mE{{\mathbb E}}

\def\mI{{\mathbb I}}

\def\mK{{\mathbb K}}

\def\mN{{\mathbb N}}

\def\mP{{\mathbb P}}

\def\mR{{\mathbb R}}

\def\sD{{\mathscr D}}

\def\sF{{\mathscr F}}

\def\sL{{\mathscr L}}

\def\sR{{\mathscr R}}
\def\sS{{\mathscr S}}

\def\E{\mathbb E}

\def\geq{\geqslant}
\def\leq{\leqslant}

\def\eps{\varepsilon}

\def\e{\mathrm{e}}

\def\p{\partial}

\def\[{{\Big[}}
\def\]{{\Big]}}
\def\<{{\langle}}
\def\>{{\rangle}}
\def\({{\Big(}}
\def\){{\Big)}}

\def\bx{{\mathbf{x}}}

\def\dif{{\mathord{{\rm d}}}}

\def\min{{\mathord{{\rm min}}}}
\def\Vol{\mathord{{\rm Vol}}}

\def\no{\nonumber}
\def\={&\!\!=\!\!&}
\def\bt{\begin{theorem}}
\def\et{\end{theorem}}
\def\bl{\begin{lemma}}
\def\el{\end{lemma}}
\def\br{\begin{remark}}
\def\er{\end{remark}}
\def\bx{\begin{Examples}}
\def\ex{\end{Examples}}
\def\bd{\begin{definition}}
\def\ed{\end{definition}}
\def\bp{\begin{proposition}}
\def\ep{\end{proposition}}
\def\bc{\begin{corollary}}
\def\ec{\end{corollary}}

\def\wt{\widetilde}

\def\geq{\geqslant}
\def\leq{\leqslant}

 \def\R{\mathbb R}
 \def\R{\mathbb R}

\def\<{\langle} \def\>{\rangle}

 \def\beq{\begin{equation}}  
 
\def\e{\text{\rm{e}}}

\allowdisplaybreaks

\begin{document}

\title[Heat kernels for non-symmetric mixed L\'evy-type operators]{Heat kernels for time-dependent non-symmetric mixed L\'evy-type operators}

\author{Zhen-Qing Chen \quad and \quad  Xicheng Zhang}

 \address{Zhen-Qing Chen:
Department of Mathematics, University of Washington, Seattle, WA 98195, USA  \\
Email: zqchen@uw.edu
 }

\address{Xicheng Zhang:
School of Mathematics and Statistics, Wuhan University,
Wuhan, Hubei 430072, P.R.China   \\
Email: XichengZhang@gmail.com
 }

\thanks{
 This work is partially supported by  Simons Foundation grant 520542  and by an NNSFC grant of China (No. 11731009). }

\begin{abstract}
In this paper we establish the existence and uniqueness of heat kernels to a large class of 
time-inhomogenous non-symmetric 
 nonlocal operators with Dini's continuous kernels. Moreover, 
  quantitative estimates including two-sided estimates, gradient estimate and fractional derivative estimate of the heat kernels
  are obtained.

\bigskip

  \noindent {{\bf AMS 2020 Mathematics Subject Classification:} 
Primary 
35K08, 
60J35, 47G20; Secondary 47D07}
	
	\medskip

  \noindent{{\bf Keywords and Phrases:} Heat kernel estimates,  
   non-symmetric nonlocal operator, Dini continuity }

\end{abstract}

\maketitle

\section{Introduction}

The purpose of this paper is to study fundamental solutions (also called heat kernels)
 and their estimates for a large class of time-inhomogenous 
non-symmetric non-local operators   on Euclidean spaces with Dini-continuous coefficients.
 
Throughout  this paper,  $\phi (r)$ is a strictly increasing continuous function on $\R_+$ with the property that 
$\phi (0)=0,\ \phi (1)=1$, and 
\begin{equation}\label{e:1.1a}
c_0^\phi := \int_0^\infty  \frac{r^2\wedge 1} {r\phi (r)} \dif r <\infty. 
\end{equation} 
Here and below, we use $:=$ as a way of definition.
Condition \eqref{e:1.1a} is necessary and sufficient for $\nu (\dif z):=\frac{1}{|z|^d \phi (|z|) }\dif z$
to be the L\'evy measure of a L\'evy process on $\R^d$.
We can divide $\phi$ into    three separate cases:
\begin{equation}\label{e:1.2}
\begin{cases}
\mbox{ Case$^\phi_1$}:   \int_{0+} \frac1{\phi (r)} \dif r <\infty,  \smallskip \\
\mbox{ Case$^\phi_2$}:  \int_{0+} \frac1{\phi (r)} \dif r =\infty \hbox{ and } \int_1^\infty \frac{1}{\phi (r)} \dif r =\infty , 
\smallskip  \\ 
\mbox{ Case$^\phi_3$}: \int_{0+} \frac1{\phi (r)} \dif r =\infty  \hbox{ and } \int_1^\infty \frac{1}{\phi (r)} \dif r <\infty,
\end{cases} 
\end{equation}
 and define
$z^{(\phi)}:\R^d \to \R^d$ by 
\begin{equation}\label{e:1.3}
z^{(\phi)}:=z \1_{\{|z|\leq 1\}}\cdot  {\1}_{\rm Case^\phi_2}+z\cdot {\1}_{\rm Case^\phi_3}.
\end{equation}
When $\phi (r)=r^\alpha$ with $ \alpha \geq 0$, \eqref{e:1.1a} holds if and only if $0<\alpha <2$; in this case,
Case$^\phi_1$, Case$^\phi_2$ and Case$^\phi_3$ correspond 
 to $ \alpha <1$, $\alpha=1$ and $\alpha >1$. 
 When $\phi (r)=r^\alpha \1_{(0, 1)} + r^\beta \1_{[1, \infty)}$ with  positive $\alpha$
and $ \beta  $,  condition \eqref{e:1.1a} holds if and only if $0<\alpha <2$; in this case,
Case$^\phi_1$ holds if and only if $\alpha <1$, Case$^\phi_2$ holds if and only if $\alpha \geq 1$ and $0<\beta \leq 1$,
and Case$^\phi_3$ holds if and only if $\alpha \geq 1$ and $\beta>1$. 

\medskip

 Let $d\geq 1$ and 
 $\kappa(t,x,z):\mR_+\times\mR^d\times\mR^d\to\mR$ be a function that is bounded between two positive constants
that is  Dini-type continuous in $x$ uniformly in $(t, z)$. 
In this paper, we study,  
 under some mild conditions on $\kappa (t, x, z)$ and $\phi (r)$, 
 the existence, uniqueness and two-sided estimates of heat kernels for  the following nonlocal 
 operator on $\R^d$:
\begin{align}\label{Non}
 \sL^\kappa_t f(x):=
   \int_{\mR^d} \left(f(x+z)-f(x)-z^{(\phi)}\cdot\nabla f(x)\right)\frac{\kappa(t,x,z)}{|z|^d \phi (|z|) }\dif z,
\end{align}
if $ \kappa (t, x, z) $ is not symmetric in  $z$, and
\begin{align}\label{Non1}
\sL^\kappa_t f(x):=\frac{1}{2}\int_{\mR^d} \left(f(x+z)+f(x-z)-2f(x)\right)\frac{\kappa(t,x,z)}{|z|^d \phi (|z|) }\dif z,
\end{align}
if $\kappa (t, x, z)$ is symmetric in $z$ for every $(t, x)\in \R_+\times \R^d$.
  For notational simplification,  
  unless otherwise specified, we define 
  \begin{align}\label{e:Delta}
\Delta^{(\phi)}_f(x,z):=
\begin{cases} \displaystyle 
\frac{f(x+z)-f(x)-z^{(\phi)}\cdot\nabla f(x)}{|z|^d \phi (|z|) }   &\hbox{if } \kappa (t, x, z) \hbox{ is not symmetric in } z, 
\medskip \\
\displaystyle 
\frac{f(x+z)+f(x-z)-2f(x)}{2|z|^d \phi (|z|) }   &\hbox{if } \kappa (t, x, z) \hbox{ is  symmetric in } z.
\end{cases}
\end{align}
Then we can write the non-local operator in \eqref{Non}-\eqref{Non1}
 in a unified way as 
\begin{equation}\label{Non2}
\sL^\kappa_t f(x) =\int_{\mR^d}\Delta^{(\phi)}_f(x,z)\, \kappa(t,x,z) \, \dif z.
\end{equation}
\begin{remark}\label{R:1.1}  \rm 
\begin{enumerate} [(i)]
\item When $\kappa (t, x, z)$ is symmetric in $z$, if $f$ is differentiable and the integral in 
\eqref{Non} is absolutely convergent, then so is the integral in \eqref{Non1} and these two integrals
give the same value.  The definition in \eqref{Non1} has the advantage
that it does not a priori require $f$ to be differentiable.

\item The reason we use $z^{(\phi)}$ in \eqref{Non} instead of the more common $z \1_{\{|z| \leq 1\}}$
in the first order correction term in \eqref{Non} is that 
this is the form 
for general $\alpha$-stable L\'evy processes where $\phi (r)=r^\alpha$ for $0<\alpha <2$
and $\kappa (t, x, z)$ is independent of $t$ and $x$.
\end{enumerate}
\end{remark} 

We consider the following conditions on $\phi$. 
There exist constants $0<c_1^\phi\leq  c_2^\phi$ , and $0<\beta_1\leq \beta_2<\infty$
 such that for all $0<r<R<\infty$,
\begin{equation}\label{e:1.1}
c^\phi_1\left(\frac{R}{r}\right)^{\beta_1} \1_{\{0<r<R\leq 1\}} \leq\frac{\phi  (R)}{\phi  (r)}\leq 
c^\phi_2\left(\frac{R}{r}\right)^{\beta_2}  .
\end{equation}
Note that the lower bound in \eqref{e:1.1} implies that for any $r\in (0, 1)$,
\begin{equation} \label{e:1.6} 
\int_r^1 \frac{1}{s\phi (s)} \dif s
 =
\frac1{\phi (r)} \int_r^1 \frac{\phi (r)}{s\phi (s)} \dif s 
\leq \frac1{c_1^\phi \phi (r)} \int_r^1 \frac{r^{\beta_1} }{s^{1+\beta_1}} \dif s 
= \frac{1-r^{\beta_1}} { c_1^\phi \beta_1 \phi (r)} 
\leq \frac1 { c_1^\phi \beta_1 \phi (r)}.
\end{equation} 
We point out that we do not assume weak lower scaling condition on $\phi$ at infinity; that is, the lower bound 
in \eqref{e:1.1} is only assumed for $0<r<R\leq 1$.
 Let 
\begin{equation}\label{e:1.21}
\gamma^{(0)}_\phi(r):=r^2\wedge1,\quad \ \gamma^{(1)}_\phi(r):=(r\wedge 1) {\1}_{\rm Case^\phi_1}
+(r^2\wedge 1){\1}_{\rm Case^\phi_2}+(r^2\wedge r ) {\1}_{\rm Case^\phi_3}.
\end{equation}
 We will also consider the condition
 \begin{align}\label{Con0}
\quad\cA^{(i)}_\phi:=
\sup_{\lambda\in(0,1]}\int^\infty_0\frac{\phi(\lambda)\gamma^{(i)}_\phi(r)}{r\phi(\lambda r)}\dif r<\infty,\quad \ i=0,1.\tag{{\bf A$^{(i)}_\phi$}}
\end{align} 
The  above condition 
is quite natural if  we want to have some approximate 
scaling properties about the heat kernel 
 (see Proposition \ref{Pr21} and Lemma \ref{Th24} below).
 Note that ({\bf A$^{(1)}_\phi$}), which will be used  in the case that $\kappa (t, x, z)$ is not symmetric in $z$, implies ({\bf A$^{(0)}_\phi$}), and ({\bf A$^{(0)}_\phi$}) implies \eqref{e:1.1a}.
It is easy to verify the following:
 \begin{enumerate} [(i)]
 \item  for  $\phi (r)=r^\alpha$,  ({\bf A$^{(1)}_\phi$}) holds
for  every $\alpha \in (0, 2)$;

\item  for $\phi (r)=r^\alpha \1_{(0, 1)} + r^\beta \1_{[1, \infty)}$, 
  ({\bf A$^{(0)}_\phi$}) holds for every $\alpha \in (0, 2)$ and $\beta >0$,
  while   ({\bf A$^{(1)}_\phi$}) holds for  every $\alpha \in (0, 2)$ and $\beta>0$ 
   except  when $\alpha=1$ and $\beta>1$. 
   \end{enumerate} 
 Some additional 
 examples are given in Section \ref{S:5} of this paper 
 that satisfy  \eqref{e:1.1} and \eqref{Con0}.
 
\medskip

We recall the following definitions about Dini and slowly varying functions.

\bd  \rm 
 Let $\ell: (0,1]\to(0,\infty)$ be a continuous function.
We call it a {\it slowly varying function}  at zero if
$$
\lim_{t\to 0}\ell(\lambda t)/\ell(t)=1 \quad \hbox{for every } \lambda>0.
$$
We call it a {\it Dini function} if $\ell$ is increasing and
$$
\int^1_0\frac{\ell(t)}{t}\dif t<\infty.
$$
We denote by $\sS_0$ (resp. $\sD_0$) the set of all slowly varying functions at zero 
 that is bounded away from zero on $[\eps,1]$ for any $\eps\in(0,1)$
 (resp.  Dini functions). For $\alpha\geq 0$, we denote by $\sR_\alpha$ 
the set of all functions $\ell(t)=t^\alpha\ell_0(t)$ for $t\in(0,1]$, where $\ell_0\in\sS_0$.
Clearly, $\sR_\alpha \subset \sD_0$ for $\alpha >0$. 
In the following, we use the convention that the definition of functions $\ell$ in $\sS_0$, $\sR_\alpha$ and $\sD_0$
are extended from $(0,1]$ to $(0,\infty)$ by setting  $\ell(t)=\ell(1)$ for $t\geq 1$. 
\ed

\bx\rm
 Clearly, $1/\ell\in\sS_0$ for $\ell\in\sS_0$, 
 and $1\in\sS_0$ but not in $\sD_0$. 
Let
$$
\ell(t):=(\log (1+ 1/(t\wedge 1)))^\alpha,\ \alpha\in\mR.
$$ 
It is easy to see that $\ell\in\sS_0$, 
and   $\ell\in\sS_0\cap\sD_0$ when   $\alpha<-1$.  
\ex

Throughout this paper,  we assume the function 
 $\kappa (t, x, z)$ in \eqref{Non} satisfies that for some $\kappa_0\geq 1$,
\begin{align}\label{Con1}
\kappa^{-1}_0\leq \kappa(t,x,z)\leq \kappa_0, \quad
 |\kappa(t,x,z)-\kappa(t,y,z)|\leq \ell^2(|x-y|),
\end{align}
where $\ell\in\sS_0\cap\sD_0$ or $\ell\in\sR_\alpha$ with $\alpha\in(0,1]$, and in the above Case$^\phi_2$,
\begin{align}\label{Con2}
\int_{  |z|\leq r} z \kappa(t,x,z)  \dif z=0  \quad \hbox{for every  }  r>0.
\end{align}
 Note that if $\kappa (t, x, z)$ is symmetric in $z$, then condition \eqref{Con2}  is automatically satisfied. 
In some situations, we will also need the following condition
\begin{align}\label{Con3}
 M^\phi_\ell(t):=\int^t_0\frac{1}{r}\left(\frac{\ell(r)}{\ell(t)}+\frac{\phi(r)}{\phi(t)}\right)\dif \phi(r)<\infty \quad \hbox{for } t\in(0,1].
\end{align} 
When      $\ell(r)=r^\eta$ and $\phi(r)=r^\alpha$ on $[0, 1]$ with  $\eta \in (0, 1]$ and $\alpha \in (0, 2)$, 
 condition  \eqref{Con3} holds if and only if 
    $\alpha>1/2$ and $\alpha+\eta>1$.
    In this case
    \begin{equation} \label{e:1.13} 
    M^\phi_\ell(t)= \left(\frac{\alpha}{\alpha+\eta-1} + \frac{\alpha}{2\alpha-1}\right)t^{\alpha-1}
    \quad \hbox{for } t\in (0, 1].
    \end{equation} 
  In comparison,  Case$^\phi_2$ and Case$^\phi_3$ correspond to $\alpha \in [1, 2)$. 
Some additional concrete  conditions for \eqref{Con3}
 to hold 
 can be found in Remark \ref{R:1.4} and in
Example \ref{E:5.6} below. 
Condition \eqref{Con3} is needed for   the gradient estimate \eqref{Grad} of the  heat kernel $p^\kappa_{t, s} (x, y)$ in 
    Theorem \ref{heat} as well as for   
 the $\sL^\kappa_t$-differentiability of the heat kernel constructed in Case$^\phi_2$ and Case$^\phi_3$
  when $\kappa (t, x,z)$ is not symmetric in $z$. 

\medskip

The study of heat kernels and their estimates is an active research area in analysis and in probability theory. 
 It has a long history for second order differential operators.
 We refer the reader to the Introduction   of \cite{Ch-Zh} for a brief history on the study of heat kernels for nonlocal operators. 
When $\phi (r)=r^\alpha$ with $\alpha \in (0, 2)$, $\kappa(t,x,z)=\kappa(x,z)$ is time-independent, symmetric in $z$
and H\"older continuous in $x$, the heat kernel of $\sL^\kappa_t$ is constructed and its 
sharp two-sided estimates,  gradient estimate and fractional derivative estimate  are obtained
in \cite{Ch-Zh}. Recently, this result has been strengthened in \cite{Ch-Zh2} by 
dropping the symmetry condition on $\kappa (x, z)$ in $z$ and allowing $\kappa $ to 
be time dependent. 
The ideas and approach of  \cite{Ch-Zh} are quite robust and they have  been adopted  to study  heat kernels for 
non-local operators  with more general L\'evy kernels;  see  \cite{KSV,  KL, GS, J,  Sz, BSK} and the references therein.
 In these works, 
$\kappa (t, x, z)$ are all assumed to be independent of $t$ and H\"older continuous in $x$, 
and,     symmetry in $z$ is assumed in  \cite{KSV,  KL, BSK}. 
  In \cite{Ch-Hu-Xi-Zh}, we studied heat kernel and its regularity and estimates for  time-inhomoegenous diffusion with jumps, 
whose infinitesimal generators have both diffusive and non-local parts. 

\smallskip

 The main feature and contributions   of this paper are
\begin{enumerate} [(i)]
\item  $\kappa (t, x, z)$ is only assumed to be Dini continuous in $x$ and can be time-inhomogenous; 

\smallskip

\item $\kappa (t, x, z)$ does not need to be symmetric in $z$; 

\smallskip
\item  the L\'evy kernel $\frac{1}{|z|^d \phi (|z|)}$ is quite general  with $\phi (r)$ satisfying \eqref{e:1.1} 
and can have light tails;  

 \smallskip
 
 \item    the lower bound and the upper bound   in our two-sided heat kernel estimate are comparable;

\smallskip

 \item We fully utilize  the rough scaling property of the non-local operator $\sL^{\kappa}_t$,   
which makes our approach in this paper more direct; 
 
\end{enumerate}

   \medskip
   
   In this paper,  we  use the following notations and conventions.  
\begin{enumerate}[$\bullet$]
\item  For $a, b\in \mR$, $a\wedge b:=\min\{a, b\}$ and $a\vee b:=\max \{a, b\}$.
Notation   $f \asymp g $ means that there are positive constants $c_1$ and $c_2$
so that $c_1 f \leq g\leq c_2 f$ on their common domain of definitions. 

\smallskip

\item The space of bounded functions on $\mR^d$ with bounded first and second
 derivatives is denoted by $C^2_b(\mR^d)$.
 
\smallskip

\item    The inverse function of $\phi (r)$ is denoted by $\phi^{-1}(r)$ . 
\smallskip

\item For $\ell\in\sD_0$, we introduce
\begin{equation}\label{e:1.14} 
\ell_\phi(t):=\ell ( \phi^{-1}(t))  \quad \hbox{and} \quad  \Gamma_{\ell}(t):=\int^t_0\frac{\ell(s)}{s}\dif s.
\end{equation} 
\item For any $T\in(0,\infty]$ and $\eps\in[0,T)$, write
$$
\mD^T_{\eps}:=\Big\{(t,x;s,y): x,y\in\mR^d \mbox{ and } s,t\geq 0 \mbox{ with } \eps<s-t<T\Big\}.
$$

\item $\Theta $ and $\Theta_1$ stand  for   sets of parameters:
$$
\Theta_1:=\left(d, \kappa_0, c^\phi_0, c^\phi_1, c^\phi_2, \beta_1, \beta_2 \right), \quad 
\Theta:=\left(\Theta_1,  \cA^{(0)}_\phi,\cA^{(1)}_\phi\right). 
$$

\item For $t>0$ and $x\in\mR^d$, we define
\begin{equation}\label{Def8}
\rho_\phi(t,x):=\frac1{ t\phi^{-1}(t)^d+|x|^d \phi (|x|) } 
\quad \hbox{and} \quad \rho_\phi(x):= \rho_\phi(1,x)= \frac1{1+|x|^d \phi (|x|) } . 
 \end{equation}
\end{enumerate}

Clearly
$$
\rho_\phi(t,x) \asymp \frac1{ t\phi^{-1}(t)^d  }  \wedge \frac1{  |x|^d \phi (|x|) }.
$$
Note that for $t\in (0, 1]$, 
\begin{equation}\label{e:1.16}
t  \int_{\mR^d}  \rho_\phi(t,x)  \dif x 
 \asymp    \int_{\{|x|\leq \phi^{-1}(t)\}} \frac1{\phi^{-1}(t)^d}  \dif x +
\int_{\{|x|> \phi^{-1}(t) \}} \frac{t}{|x|^d \phi (|x|)}  \dif x    
 \asymp  1+ \int_{\phi^{-1}(t)}^\infty \frac{t}{r\phi (r) } \dif r .
 \end{equation} 
 By \eqref{e:1.1a} and \eqref{e:1.6},   for  $t \in (0, 1]$,
 $$
  \int_{\phi^{-1}(t)}^\infty \frac{t}{r\phi (r) } \dif r 
 =\int_1^\infty \frac{t}{r\phi (r) } \dif r + \int_{\phi^{-1}(t)}^1 \frac{t}{r\phi (r) } \dif r \lesssim 1 . 
 $$
 Thus  we have by \eqref{e:1.16},
  \begin{equation}\label{e:1.17} 
    \int_{\mR^d}  \rho_\phi(t,x)  \dif x  \asymp 1/t  \quad \hbox{for } t\in (0, 1].
    \end{equation}
Since $\phi$ is increasing, 
 $$
t\phi^{-1}(t)^d+|x|^d \phi (|x|)\leq 2\left(\phi^{-1}(t)+|x|\right)^{d}\left(\phi \big(\phi^{-1}(t)+|x|\big)\right),
$$
and since $\phi(2r)\leq c\phi(r)$,
$$
\left(\phi^{-1}(t)+|x|\right)^{d}\left(\phi \big(\phi^{-1}(t)+|x|\big)\right)\leq c \left(t\phi^{-1}(t)^d+|x|^d \phi (|x|)\right).
$$
Thus 
\begin{align}\label{EK0}
\rho_\phi(t,x) \asymp \frac1{\left(\phi^{-1}(t)+|x|\right)^d  \phi \big(\phi^{-1}(t)+|x| \big)}.
\end{align}
It is known from \cite{CK08}  that the transition density function $p(t, x, y)$
for a pure jump symmetric L\'evy process on $\R^d$ with L\'evy measure $\frac1{|z|^d \phi (|z|)} \dif z$  has
the two-sided estimates:  
$$
p (t, x, y)\asymp   t\rho_{\phi} (t, |x-y|)\quad \hbox{for all } t>0 \hbox{ and } x, y\in \R^d.
$$
We will show in this paper  
that the above estimate also holds for purely discontinuous
  non-symmetric L\'evy processes whose L\'evy measure is comparable to that of isotropic L\'evy process on $\R^d$.
  In fact, more is true; see Theorem \ref{Th1} for a precise statement.

The following is the main result of this paper.

\bt\label{heat}
Suppose that one of the following two assumptions holds:
\begin{enumerate}[{\bf (H1)}]
\item If $\kappa(t,x,z)=\kappa(t,x,-z)$, we assume 
 \eqref{e:1.1}, ({\bf A$^{(0)}_\phi$}) and \eqref{Con1}.
\item If $\kappa(t,x,z)\not=\kappa(t,x,-z)$, we assume 
 \eqref{e:1.1}, ({\bf A$^{(1)}_\phi$}), \eqref{Con1}
 and \eqref{Con2}.
   In addition, we assume \eqref{Con3} for Case$^\phi_2$ and Case$^\phi_3$. 
   \end{enumerate}
 Then there is a unique continuous function $p^\kappa_{t,s}(x,y)=p^\kappa(t,x;s,y)$ on $\mD^\infty_0$
(called the fundamental solution or heat kernel of $\sL^\kappa_t$) 
so that  for all $f\in C^2_b(\mR^d)$, 
$$
P^\kappa_{t,s}f(x):=\int_{\mR^d}p^\kappa_{t,s}(x,y)f(y)\dif y
$$
has the following properties: $\sL^\kappa_r P^\kappa_{t,s}f(x)$ exists for each $s>t\geq 0$ and $x\in \R^d$,
\begin{align}\label{EQ}
P^\kappa_{t,s}f(x)=f(x)+\int^s_t\!\sL^\kappa_r P^\kappa_{r,s}f(x)\dif r \quad \hbox{for every }  (t,x)\in[0,s]\times\mR^d,
\end{align}
 and
\begin{enumerate}[\rm (i)]

\item 
 for each $t_0\in[0,s)$ and $x\in\mR^d$, 
it holds that
$$
\lim_{t\downarrow t_0}\left|\sL^\kappa_tP^\kappa_{t,s}f(x)-\sL^\kappa_tP^\kappa_{t_0,s}f(x)\right|=0 ;
$$
 
 \item when $\kappa (t, x, z)$ is not symmetric in $z$,   
  $x\mapsto\nabla P^\kappa_{t,s}f(x)$   exists and is 
   continuous on $\mR^d$ in Case$^\phi_2$ and Case$^\phi_3$ for each $0\leq t<s$;

\item 
for any bounded and uniformly continuous function $f$ on $\R^d$, 
      \begin{align}
        \lim_{|t-s|\to 0}\|P^\kappa_{t,s}f-f\|_\infty=0.\label{cz1}
      \end{align}

\end{enumerate}
Moreover, the heat kernel $p^\kappa_{t,s}(x,y)$ enjoys the following properties:
 \begin{enumerate}[\bf (a)]
\item (Two-sided estimate) For any $T>0$, there is a 
$c_1=c_1(T,\ell,\Theta)>1$,
 such that on $\mD^T_0$,
\begin{align}\label{GR1}
c_1^{-1}(s-t)\rho_\phi (s-t,x-y)\leq p^{\kappa}_{t,s}(x,y)\leq c_1(s-t)\rho_\phi(s-t,x-y).
\end{align}

\item (Fractional derivative estimate) 
For any $T>0$, there is a constant 
 $c_2=    c_2(T,\ell,\Theta) > 0 $ such that  on $\mD^T_0$,
\begin{align}\label{Fra} 
 \int_{\mR^d}\big|\Delta^{(\phi)}_{p^\kappa_{t,s}(\cdot,y)}(x,z)\big|
 \,  \dif z 
 \leq c_2\left(\frac{\Gamma_{\ell_\phi}(s-t)}{\ell_\phi(s-t)}\right)
\rho_\phi(s-t,x-y). 
\end{align}

\item (Gradient estimate) 
Suppose that \eqref{Con3} holds. 
 Then $x\mapsto  p^\kappa_{t,s}(x,y) $ is continuously differentiable for each $0\leq t<s$, and
for every $T>0$ there is a constant 
$c_3 = c_3 (T,\ell,\Theta)  >0$ so that on $\mD^T_0$,
\begin{align}\label{Grad}
 |\nabla p^\kappa_{t,s}(\cdot,y)(x)|\leq 
 c_3\left(\frac{s-t}{\phi^{-1}(s-t)}+M^\phi_\ell\circ\phi^{-1}(s-t)\right)
\rho_\phi(s-t,x-y).  
\end{align}
 
\item
(Conservativeness)
For every $0<t<s$ and $x\in \mR^d$, 
$$
\int_{\mR^d}p^\kappa_{t,s} (x,y)\dif y=1. 
$$

\item (Chapman-Kolmogorov equation) 
For all $0<t<r<s<\infty$ and $x,y\in\mR^d$, 
\begin{align}\label{eq21}
\int_{\mR^d}p^\kappa_{t,r} (x,z)p^\kappa_{r,s} (z,y)\dif z=p^\kappa_{t,s} (x,y).
\end{align}

\item (Generator) 
For any $f\in C_b^2(\mR^d)$, we have
      \begin{align}
        P^\kappa_{t,s}f(x)=f(x)+\int^s_t\!P^\kappa_{t,r}\sL^\kappa_rf(x)\dif r.\label{eqge}
      \end{align}
\end{enumerate}
\et

 \br \label{R:1.4} \rm 
\begin{enumerate}[(i)]
   \item If $\ell\in\sR_\alpha$ for some $\alpha>0$, then 
$\lim_{t\to 0}  {\Gamma_{\ell_\phi}(t)}/{\ell_\phi(t)} = 1/{\alpha}$ by \eqref{Ka} below.
In this case,  \eqref{Fra}   implies that 
\begin{align}\label{Frab}
 |\sL^\kappa_t p^\kappa_{t,s}(\cdot,y)(x)|\leq \wt c_2 
\rho_\phi(s-t,x-y)\quad  \hbox{on } \mD^T_0.
\end{align}

\item We emphasize that the jumping kernel $j(t, x, z):=\frac{\kappa(t,x,z)}{|z|^d \phi (|z|) }$
for $\sL^{\kappa}_t$  of  \eqref{Non}  can have light tail in $z$ at infinity.
For instance,  in the example when $\phi (r)=r^\alpha \1_{(0, 1)} + r^\beta \1_{[1, \infty)}$,
one can check (see Example \ref{E:5.3} below)  that \eqref{e:1.1} 
holds with $\beta_1=\alpha$ and $\beta_2=\alpha\vee\beta$ and 
 ({\bf A$^{(0)}_\phi$}) 
is  satisfied 
for any $\alpha\in(0,2)$ and  $\beta>0$; while 
 ({\bf A$^{(1)}_\phi$}) 
holds for all   $\alpha\in(0, 2)$ and  $\beta>0$ except when $\alpha =1$ and $\beta >1$.
In other words, $\beta$ can be any positive number and so   $\beta$ 
  can be larger than or equal 2. 
Two-sided heat kernel estimates  for symmetric pure jump processes with light polynomial decay jumping kernels
have recently been studied in  \cite{BKKL, CKW}. 
In particular, see some symmetric analogous estimates of \eqref{GR1} in \cite[(1.4)]{CKW} as
well as \cite[Theorem 1.2 and Theorem 1.4(i)]{BKKL}.

  \item 
 As mentioned earlier, when 
  $\ell(s)=s^\beta$ and $\phi(s)=s^\alpha$ on $[0, 1]$ for $\beta \in (0, 1]$ and $\alpha \in (0, 2)$,
   condition \eqref{Con3} holds if and only if $\alpha>1/2$ and 
   $\alpha+\beta>1$.
In this case,  $M^\phi_\ell(t)$ is given by \eqref{e:1.13}.
So the gradient estimate \eqref{Grad} takes the following form
  \begin{align}\label{Grada}
 |\nabla p^\kappa_{t,s}(\cdot,y)(x)|\leq 
   \wt   c_3  (s-t)^{1-(1/\alpha)}  
\rho_\phi(s-t,x-y)\quad  \hbox{on } \mD^T_0,
\end{align}
  which recovers and extends the gradient estimate in \cite[Theorem 1.1(5)]{Ch-Zh} and \cite[Theorem 1.1(v)]{Ch-Zh2}.
   Gradient estimate for $\phi (r)=r^\alpha$ with $\alpha \leq 1/2$ would need more restrictive assumption on the kernel $\kappa (t, x, z)$; 
 see \cite{LSX}.
 Condition \eqref{Con3} is also satisfied when $\ell \in \sD_0 \cap \sS_0$ and $\phi (r)=r$.
 In this case, 
\begin{equation}\label{e:1.31} 
 M^\phi_\ell (t) =\frac{1}{\ell (t)} \int_0^t \frac{\ell (r)}r \dif r + 1 \asymp \frac{\Gamma_\ell (t)}{\ell (t)}
 \end{equation} 
 in view of Proposition \ref{Pr32}(ii) and so    the gradient estimate \eqref{Grad} has the form   
  \begin{align}\label{Gradc}
 |\nabla p^\kappa_{t,s}(\cdot,y)(x)|\leq 
   \wt   c_3   \frac{\Gamma_\ell (s-t)}{\ell (s-t)}  
\rho_\phi(s-t,x-y)\quad  \hbox{on } \mD^T_0.
\end{align}
 
  \item  Recall $\beta_1$ is the exponent
  in \eqref{e:1.1} for $\phi$.  
  More generally,  when $\ell \in \sR_\eta$ for some $\eta\geq 0$, we will show in Example \ref{E:5.6}  below that 
    condition \eqref{Con3} holds if $\beta_1 >1/2$ and $\beta_1 + \eta >1$;
   in this case, $M^\phi_\ell (t) \asymp \phi (t)/t$ on $(0, 1]$, and consequently 
   the gradient estimate \eqref{Grad}   takes the following form: 
  \begin{align}\label{Gradb}
 |\nabla p^\kappa_{t,s}(\cdot,y)(x)|\leq 
  \wt   c_3 \frac{  s-t  } {\phi^{-1}(s-t)}     
\rho_\phi(s-t,x-y)\quad \hbox{on } \mD^T_0. 
\end{align}
This extends the gradient estimate \cite[Theorem 1.2(4)]{KSV}, where $\frac{1}{|z|^d \phi (|z|)} \dif z$ is assumed to be
the L\'evy measure of a subordinate Brownian motion, 
the lower scaling exponent $\beta_1$ of $\phi$ in \eqref{e:1.1} is within $(2/3, 2)$, 
  and $\kappa (t, x, z)$ is time-independent, symmetric in $z$ and uniformly $\eta$-H\"older continuous in $x$. 
   with $\beta_1 + \eta >1$,
   
 \item Our approach exploits  the rough scaling property of the operator $\sL^{\kappa}_t$; see \eqref{e:2.4}
 and Proposition \ref{Pr21}.  This allows us to reduce the study of heat kernel $p^\kappa_{t,s}(x, y)$ for  general 
 time $0\leq t<s\leq T$  to $p^\kappa_{0, 1}(x, y)$. In particular, this rough scaling property combined  with 
 an idea from T. Watanabe \cite{Wa} allows us to derive 
  two-sided estimates 
 as well as derivative estimates 
 for heat kernel estimate of time-dependent 
 L\'evy process in Theorem \ref{Th1} via \eqref{Sca}, which is of independent interest and plays a key role in our investigation
 of heat kernels of the space and time dependent non-local operator $\sL^\kappa_{t}$.
 
 \item Under conditions {\bf (H1)} and \eqref{Con3},   one can see from the proof of Theorem \ref{heat} 
 (in particular, Lemma \ref{Th24})  below that in fact 
 \eqref{Fra} also holds for $\Delta^{(\phi)}_{p^\kappa_{t,s}}$ being defined using the first expression in  \eqref{e:Delta}.
 Consequently, $x\mapsto p^\kappa_{t,s} (x, y)$ is pointwisely $\sL^\kappa_t$-differentiable  for every fixed $s>t\geq 0$
 and $y\in \mR^d$ with $\sL^\kappa_t$ being defined by \eqref{Non}, and Theorem \ref{heat}(i) holds in this sense as well.
 \end{enumerate}
\er
 
   \medskip
   
 The rest of the paper is organized as follows. 
 In Section \ref{S:2}, we study   heat kernel estimates for $\sL^\kappa$ when $\kappa$ does not depend on
 the state variable $x$,
 or equivalently, transition density functions of time-inhomogeous L\'evy processes. 
  In particular, the derivative estimates as well as the continuous dependence of the heat kernel in $\kappa$ are derived.
 In Section  \ref{S:3}, the properties of slowly varying functions and the basic convolution inequality are presented.
 In Section \ref{S:4}, we prove our main result Theorem \ref{heat} 
 using the classical Levi method   in time-inhomogeous and non-local operator setting.
 In Section \ref{S:5}, several examples are provided 
to illustrate the main result of this paper and its scope. 
  
\section{Heat kernel estimates of $\sL^\kappa_t$ with $\kappa(t,x,z)=\kappa(t,z)$}\label{S:2} 

Throughout this section, 
$\phi$ is a strictly increasing continuous function on $\mR_+$ with $\phi (0)=0$ and $\phi (1)=1$
satisfying conditions \eqref{e:1.1a} and \eqref{e:1.1}, and 
$\kappa (t, x, z)=\kappa (t, z)$ is  independent of $x$ and 
we assume that for some $\kappa_0\geq 1$ and all $r>0$,
\begin{equation}\label{Ass1a}
\kappa^{-1}_0\leq \kappa(t,z)\leq \kappa_0,\quad  {\1}_{\rm Case^\phi_2}\int_{|z|\leq r}z\kappa(t,z)\dif z=0.
\end{equation}
Note that if $\kappa (t, z)$ is symmetric in $z$, then $\int_{|z|\leq r}z\kappa(t,z)\dif z=0$ is automatically satisfied.

\subsection{Scaling property}
 Let $N(\dif t,\dif z)$ be a time-inhomogenous Poisson random measure  
with intensity measure $\frac{\kappa(t,z)}{|z|^d \phi (|z|) }\dif z\dif t$.
Define
\begin{equation}\label{e:2.3a}
X^\kappa_{t,s}:=\int^s_t\!\!\!\int_{\mR^d}z \tilde N(\dif r,\dif z)+\int^s_t\!\!\!\int_{\mR^d}(z-z^{(\phi)})\frac{\kappa(r,z)}{|z|^d \phi (|z|) }\dif z\dif r,
\end{equation}
where $\tilde N(\dif t,\dif z):=N(\dif t,\dif z)-\frac{\kappa(t,z)}{|z|^d \phi (|z|) }\dif z\dif t$.
Note that the process $X^\kappa_{t,s}$ is a time-inhomogenous L\'evy process on $\mR^d$
 in the sense that it has independent increments.

Sometimes we write $X^{\kappa, \phi}_{t,s}$ for $X^\kappa_{t,s}$ if we want to emphasize its dependence on $\phi$ as well. 
By It\^o's formula, we have
$$
\mE f(X^\kappa_{t,s})=\mE\int^s_t\sL^{\kappa}_r f(X^\kappa_{t,r})\dif r,\quad f\in C^2_b(\mR^d).
$$
In particular, if we take $f(x)=\e^{\mathrm{i}\xi\cdot x}$, then one finds that  the characteristic function of $X^\kappa_{t,s}$ is given by
$$
\mE \e^{\mathrm{i}\xi\cdot X^\kappa_{t,s}}=\exp\left( \int^s_t\!\!\!\int_{\mR^d}(\e^{\mathrm{i}\xi\cdot z}-1-\mathrm{i}\xi\cdot z^{(\phi)})
\frac{\kappa(r,z)}{|z|^d \phi (|z|) }\dif z\dif r\right) .
$$
By the definition of $z^{(\phi)}$ and a change of variable, as well as \eqref{Ass1a},
we conclude from the last display that for every $\lambda >0$,
\begin{equation}\label{e:2.4}
\left\{ (\phi^{-1}(\lambda))^{-1 } X^{\kappa, \phi}_{\lambda t,  \lambda s} , s>t \right\} \ 
\hbox{ has the same distribution as } \
\Big\{ X^{\kappa_\lambda ,  \phi_\lambda}_{t, s}, s>t \Big\},
\end{equation}
 where 
 $$
\kappa_\lambda (r, z)
:= \kappa (\lambda r, \phi^{-1}(\lambda) z) \quad \hbox{ and } \quad 
\phi_\lambda ( r):= \phi (\phi^{-1}(\lambda) r)/ \lambda,
 $$ 
and  
\begin{align}\label{Ch}
\mE \e^{\mathrm{i}\xi\cdot X^\kappa_{t,s}}
=\exp\left((s-t)\int^1_0\!\!\!\int_{\mR^d}(\e^{\mathrm{i}\xi\cdot z}-1-\mathrm{i}\xi\cdot z^{(\phi)})
\frac{\kappa(t+(s-t)r,z)}{|z|^d \phi (|z|) }\dif z\dif r\right).
\end{align}
By  \eqref{Ass1a} and \eqref{e:1.1},   for $|\xi|\geq 1$, letting $\bar\xi:=\xi/|\xi|$, we have
\begin{align}\label{DH1}
|\mE \e^{\mathrm{i}\xi\cdot X^\kappa_{t,s}}|
 &\leq   \exp\left( -   \frac{ s-t}{  \kappa_0}   \int_{ \mR^d } 
     \frac{1-\cos (\xi\cdot z)  }{|z|^{d}\phi(|z|)}\dif z \right) = \exp\left( -   \frac{s-t}{  \kappa_0}   \int_{ \mR^d } 
     \frac{1-\cos (\bar\xi\cdot z)  }{|z|^{d}\phi(|z|/|\xi|)}\dif z \right)\no\\
&\leq   \exp\left( -   \frac{c_1^\phi|\xi|^{\beta_1}(s-t)}{  \kappa_0}   \int_{|z|\leq 1} 
     \frac{1-\cos (\bar\xi\cdot z)  }{|z|^{d}\phi(|z|)}\dif z \right)\leq   \exp\left( - c|\xi|^{\beta_1}(s-t)\right), 
\end{align}
  where $c=c(\Theta_1)>0$.
Hence, $X^{\kappa,\phi}_{t,s}$ admits a smooth   density function $p^{\kappa,\phi}_{t,s}(x)$ given by 
the inverse Fourier transform 
\begin{align}\label{Den1}
p^{\kappa,\phi}_{t,s}(x )
 =(2\pi)^d\int_{\mR^d}\e^{-\mathrm{i}x\cdot\xi}\mE \e^{\mathrm{i}\xi\cdot X^{\kappa,\phi}_{t,s}}\dif\xi
=(2\pi)^d\int_{\mR^d}\mE \e^{\mathrm{i}\xi\cdot (X^{\kappa,\phi}_{t,s}-x)}\dif\xi.
\end{align}
 Moreover, with $\bar p^{\kappa,\phi}_{t,s}(x):=p^{\kappa,\phi}_{t,s}(-x)$, 
\begin{align}\label{DL6}
\p_t\bar p^{\kappa,\phi}_{t,s}(x)+\sL^\kappa_t\bar p^{\kappa,\phi}_{t,s}(x)=0 \ \hbox{ for }  s>t 
\  \hbox{ with } \ \lim_{t\uparrow s}\bar p^{\kappa,\phi}_{t,s}(x) \dif x =\delta_0 (\dif x),
\end{align}
where the limit is taken in the weak sense. 

By  \eqref{e:2.4}, we have the following scaling property, which will play a basic role in the sequel and simplify many calculations.

\bp\label{Pr21} 
 For $0\leq t<s\leq 1$, define
$$
\wt\kappa(r,z):=\kappa(t+(s-t)r, \phi^{-1}(s-t)z),\quad \wt \phi (u) := \phi (u \phi^{-1}(s-t))/ (s-t).
$$ 
 
\begin{enumerate}[\rm (i)]
\item $\wt\kappa$ satisfies \eqref{Ass1a} with the same $\kappa_0$,
 and $\wt \phi$ satisfies \eqref{e:1.1} with the same parameters as $\phi$.
 
\item For $x\in \mR^d$, it holds
\begin{align}\label{Sca}
p^{\kappa, \phi}_{t,s}(x)= \left( \phi^{-1} (s-t)\right)^{-d }p^{\wt\kappa, \wt \phi}_{0,1}( x/\phi^{-1}(s-t)),
\end{align}
and for another bounded measurable 
$\kappa'$  satisfying \eqref{Ass1a},  
$$
\left(\sL^{\kappa',\phi}p^{\kappa, \phi}_{t,s}\right)(x)=(s-t)^{-1} \left( \phi^{-1} (s-t)\right)^{-d }
\left(\sL^{\wt\kappa',\wt\phi}p^{\wt\kappa, \wt \phi}_{0,1}\right)( x/\phi^{-1}(s-t)),
$$
where $\wt\kappa'(z):=\kappa'(\phi^{-1}(s-t)z)$.

\end{enumerate}
\ep

\smallskip

\subsection{Two-sided estimate of $p^{\kappa,\phi}_{0,1}$}
In this subsection we show the sharp two-sided estimate of $p^{\kappa,\phi}_{0,1}$ by a purely probabilistic argument.
Recall that the function $\rho_\phi (x)$ is defined by \eqref{Def8}.

\bt\label{Th1}
Under \eqref{Ass1a}, there is a constant
$c_0=c_0(\Theta_1)>1$ such that for all $x\in\mR^d$,
\begin{align}\label{ER11}
c^{-1}_0 \rho_\phi(x)\leq p^{\kappa,\phi}_{0,1}(x)\leq c_0\ \rho_\phi (x),
\end{align}
and for each $j\in\mN$, there is a constant $c_j=c_j(\Theta_1)>0$ 
so that for all $x\in \mR^d$, 
 \begin{align}\label{ER101}
|\nabla^jp^{\kappa,\phi}_{0,1}(x)|\leq c_j\, \rho_\phi (x).
\end{align}
\et

\medskip

Define
\begin{align*}
&\Upsilon_1 
:=\int^1_0\!\!\!\int_{|z|\leq 1}z \wt N(\dif r,\dif z)+\int^1_0\!\!\!\int_{|z|\leq 1}(z-z^{(\phi)})\frac{\kappa(r,z)}{|z|^d \phi (|z|) }\dif z\dif r,\\
&\Upsilon_2 
:=\int^1_0\!\!\!\int_{|z|>1}z \wt N(\dif r,\dif z)+\int^1_0\!\!\!\int_{|z|>1}(z-z^{(\phi)})\frac{\kappa(r,z)}{|z|^d \phi (|z|) }\dif z\dif r.
\end{align*}
 Note that $\Upsilon_1$ and $\Upsilon_2$
are independent and have the characteristic functions
\begin{align}
\mE \e^{\mathrm{i}\xi\cdot \Upsilon_1}&=\exp\left(\int^1_0\!\!\!\int_{|z|\leq 1}(\e^{\mathrm{i}\xi\cdot z}-1-\mathrm{i}\xi\cdot z^{(\phi)})
\frac{\kappa(r,z)}{|z|^d \phi (|z|) }\dif z\dif r\right) =:\e^{\varphi_1(\xi)},\label{Den}\\
\mE \e^{\mathrm{i}\xi\cdot \Upsilon_2}&=\exp\left( \int^1_0\!\!\!\int_{|z|>1}(\e^{\mathrm{i}\xi\cdot z}-1-\mathrm{i}\xi\cdot z^{(\phi)})
\frac{\kappa(r,z)}{|z|^d \phi (|z|) }\dif z\dif r\right) =:\e^{\varphi_2(\xi)}.\label{Den0}
\end{align}
By a similar calculation as that for \eqref{DH1}, we have
$$ 
\left| \mE \e^{\mathrm{i}\xi\cdot \Upsilon_1} \right| = \left| \e^{\varphi_1 (\xi)} \right| \leq \e^{-c  | \xi|^{\beta_1}}.
$$
Consequently,
$\Upsilon_1$ has a smooth density $p_1 (x)$ with respect to the Lebesgue measure on
$\R^d$. Since $X^\kappa_{0, 1}$ is the independent sum of $\Upsilon_1$ and $\Upsilon_2$, we have
\begin{align}\label{ER5}
p^{\kappa,\phi}_{0,1}(x)=\mE \left[ p_1 (x-\Upsilon_2) \right].
\end{align}
To get the two-sided estimate of $p^\kappa_1(x)$, we prepare the following two lemmas.

\bl\label{Le22}
{\rm (i)} For any $R>0$, there is a $\delta =\delta (R,  \Theta_1)>0 $  so that 
$$
\inf_{x\in B_R}p_1 (x)\geq\delta.
$$
{\rm (ii)} For any integer $m,j\in\mN_0$, there is a constant 
 $c=c( m, j, \Theta_1)>0$ 
so that 
$$
|\nabla^jp_1(x)|\leq c(1+|x|^2)^{-m}.
$$
\el

\begin{proof}
(i) Let $\Upsilon_{11}$ and $\Upsilon_{12}$ be two independent random variables with the characteristic functions
\begin{align}
\mE \e^{\mathrm{i}\xi\cdot \Upsilon_{11}}&=\exp\left( \int_{|z|\leq 1}
\left(\e^{\mathrm{i}\xi\cdot z}-1-\mathrm{i}\xi\cdot z^{(\phi)} \right)
\left( \frac{\kappa(z) }{|z|^d \phi (|z|) } - \frac{\kappa_0c^\phi_1}{2 |z|^{d+\beta_1} }\right) \dif z\right)
=:\e^{\varphi_{11}(\xi)},\label{ER1}\\
\mE \e^{\mathrm{i}\xi\cdot \Upsilon_{12}}&=\exp\left(\int_{|z|\leq 1}(\e^{\mathrm{i}\xi\cdot z}-1-\mathrm{i}\xi\cdot z^{(\phi)})
\frac{\kappa_0c^\phi_1}{2 |z|^{d+\beta_1} }\dif z\right)=:\e^{\varphi_{12}(\xi)},\no
\end{align}
where $\kappa(z):=\int^1_0\kappa(r,z)\dif r\geq\kappa_0$ by \eqref{Ass1a}.
Let $p_{11}$ and $p_{12}$ be the continuous density functions of $\Upsilon_{11}$ and $\Upsilon_{12}$, respectively.
 Clearly,  
\begin{align}\label{ER92}
p_1(x)=\int_{\mR^d}p_{11}(x-z)p_{12}(z)\dif z.
\end{align}
Since $\Upsilon_{12}$ is a truncated rotationally symmetric $\beta_1$-stable random variable, it is well known
(see, e.g., \cite{CKK})  that $p_{12}$ is strictly positive
on $\mR^d$. On the other hand, we have by \eqref{ER1} and \eqref{DH1}  that 
$$
\mE|\Upsilon_{11}|\leq  c(\Theta_1) <\infty.
$$
Hence, by \eqref{ER92}, we have for any $R_1>R$ and $x\in B_{R}$, 
\begin{align*}
p_1(x)&=\int_{\mR^d}p_{11}(x-z)p_{12}(z)\dif z\geq\inf_{z\in B_{{R_1}}}p_{12}(z)\int_{|z|\leq {R_1}}p_{11}(x-z)\dif z\\
&=\inf_{z\in B_{{R_1}}}p_{12}(z)\Big(1-\mP(|\Upsilon_{11}-x|>{R_1})\Big)\\
&\geq\inf_{z\in B_{{R_1}}}p_{12}(z)\Big(1-(\mE|\Upsilon_{11}| +R)/{R_1}\Big),
\end{align*}
which yields (i) by taking ${R_1}$ large enough.

\medskip

(ii) Using the inverse Fourier    transform, 
for every integer $m\geq 1$,  we have by \eqref{Den}
\begin{eqnarray*}
(1+|x|^2)^m \,  |  \nabla^j p_1 (x) | 
\leq  (2\pi)^d \int_{\R^d} |\xi|^j \left|  (\mI-\Delta)^m \e^{\varphi_1(\xi)} \right| \dif \xi\leq 
 c  (\Theta_1) <\infty. 
\end{eqnarray*}
The proof is complete.
\end{proof}
 
\bl\label{Le23}
For any $R>2$, there is a 
 constant $c_1=c_1(R, \Theta_1)>0$ so that
 for all $x\in\mR^d$,
\begin{align}\label{EQ1}
\frac{c_1^{-1}}{ (1+|x|)^{d} \phi (1+|x|)} \leq \mP(\Upsilon_2\in B_R(x))\leq \frac{c_1}{ (1+|x|)^{d} \phi (1+|x|)}  .
\end{align}
\el

\begin{proof}
Observe that by \eqref{Den0},
\begin{align*}
\mE \e^{\mathrm{i}\xi\cdot \Upsilon_2}&=\exp\left(\int_{\mR^d}(\e^{\mathrm{i}\xi\cdot z}-1)\nu(\dif z)\right)
\exp\left( -\mathrm{i}\xi\cdot b\right),
\end{align*}
where $\nu(\dif z):=   \1_{\{|z|>1\}} \left(\int^1_0\kappa(r,z) \dif r\right)     \frac{1}{|z|^d \phi (|z|)} \, \dif z$ and $b:=\int_{\mR^d}z^{(\phi)}\nu(\dif z)$.
Let $\eta:=\{\eta_n,n\in\mN\}$ be a family of i.i.d. random variables in $\mR^d$ with distribution $\nu/\lambda$,  
where 
$$
\lambda:=\nu(\mR^d)\leq c(\Theta_1 )<\infty.
$$ 
Let $S_0=0$ and 
$S_n:=\eta_1+\cdots+\eta_n.$
Let $N$ be a Poisson random variable with parameter $\lambda$, which is independent of $\eta$.
It is easy to see that
$$
S_N\stackrel{(d)}{=}\Upsilon_2+b.
$$
Now, by definition we have
\begin{align*}
\mP(\Upsilon_2\in B_R(x))&=\mP(S_N\in B_R(x+b))
=\sum_{n=1}^\infty\mP\Big(S_n\in B_R(x+b)\Big)\mP(N=n)\\
&=\e^{-\lambda}\sum_{n=1}^\infty\frac{1}{n!}\int_{\mR^{nd}}\1_{\sum_{j=1}^n z_j\in B_R(x+b)}\nu(\dif z_1)\cdots\nu(\dif z_n).
\end{align*} 
When $|x+b|<R+1$, 
the upper bound in \eqref{EQ1} for 
$ \mP(\Upsilon_2\in B_R(x))$ trivially holds. 
Thus we assume that $|x+b|\geq R+1$.
Notice that $\sum_{j=1}^n z_j\in B_R(x+b)$ implies that there is at least one $i$ such that $|z_i|>(|x+b|-R)/n$. Hence,
$$
\mP(\Upsilon_2\in B_R(x))\leq\e^{-\lambda}\sum_{n=1}^\infty\frac{1}{n!}
\left(\sum_{i=1}^n \int_{\mR^{nd}}\1_{\sum_{j=1}^n z_j\in B_R(x+b)}\1_{|z_i|>(|x+b|-R)/n}\nu(\dif z_1)\cdots\nu(\dif z_n)\right).
$$
Recalling $\nu(\dif z_i)=\1_{|z_i|>1}(|z_i|^d\phi(|z_i|))^{-1}\left(\int^1_0\kappa(r,z_i)\dif r\right)\dif z_i$, we have by   \eqref{Ass1a}
\begin{align*}
\mP(\Upsilon_2\in B_R(x))
&\leq \e^{-\lambda}  \sum_{n=1}^\infty
 \frac{ \kappa_0 n^d }{ (|x+b|-R)^{d} {\phi ((|x+b|-R)/n)}}  \frac{1}{n!}\\
&\quad\times\left(\sum_{i=1}^n\int_{\mR^{nd}}\1_{\sum_{j=1}^n z_j\in B_R(x+b)}\nu(\dif z_1)\cdots\dif z_i\cdots\nu(\dif z_n)\right)\\
&\leq \e^{-\lambda} \sum_{n=1}^\infty
\frac{ \kappa_0 c^\phi_2 n^{d+\beta_2}  }{ (|x+b|-R)^{d} \phi ( |x+b|-R )}   |B_R|
\frac{\lambda^{n-1}}{n!} \\ 
&\leq \frac{c_1  (\Theta_1 )  }{ (|x+b|-R)^{d} \phi ( |x+b|-R )}  \, |B_R|
\leq \frac{c_2  (R,\Theta_1 )}{(1+|x|)^{d} \phi (1+ |x| )},
\end{align*}
where in the second inequality we used \eqref{e:1.1} and the translation invariance property of the Lebesgue measure in $z_i$-variable.
On the other hand, for any $x\in \R^d$, since $R>2$,
\begin{align*}
&\mP(\Upsilon_2\in B_R(x))\geq\e^{-\lambda}\int_{\mR^{d}}\1_{\{z_1\in B_R(x+b)\}}\nu(\dif z_1)\\
&\quad\geq \frac{ \kappa_0^{-1}\e^{-\lambda}|B_R(x+b)\cap B^c_1| }{(|x|+|b|+R)^d \phi ((|x|+|b|+R) }
 \geq  \frac{c_2  (R,\Theta_1 )}{(1+|x|)^{d} \phi (1+ |x| )}.
 \end{align*}
Combining the above calculations, we get the desired estimate.
\end{proof}

Now we can give

\begin{proof}[Proof of Theorem \ref{Th1}]
Our proof is adapted from \cite{Wa}. Let $R>2$.
For the lower bound, by (i) of Lemma \ref{Le22}, we have 
$$
\delta:=\inf_{z\in B_R}p_1 (z)>0.
$$
Hence, by \eqref{ER5} and Lemma \ref{Le23},  
\begin{equation}\label{e:2.17}
p^{\kappa,\phi}_{0,1}(x)=\E \left[ p_1(x-\Upsilon_2)\right]\geq \delta\mP(|x-\Upsilon_2|\leq R) \geq \frac{\delta c^{-1}_1}{ (1+|x|)^d \phi (1+|x|) }  .
\end{equation}
For the upper bound, by \eqref{ER5} again, we have
 \begin{align}\label{ER6}
p^{\kappa,\phi}_{0,1}(x)\leq\mE \Big[ p_1(x-\Upsilon_2)\1_{|x-\Upsilon_2|\leq|x|/2}\Big] +\sup_{|z|>|x|/2}p_1(z).
\end{align}
By (ii) of Lemma \ref{Le22}, we can choose $N$-points $z_1,\cdots,z_N\in B_{|x|/2}$ and $\eps>0$ such that 
$$
 B_{|x|/2}\subset \cup_{j=1}^N B_\eps(z_j) \quad \hbox{ and } \quad \sum_{j=1}^N\sup_{z\in B_\eps(z_j)}p_1 (z)\leq c_4,
$$
where $c_4$ only depends on $\eps,\kappa_0,d,\alpha$.
Hence, by Lemma \ref{Le23}, we have
\begin{align*}
\mE \Big[ p_1 (x-\Upsilon_2)\1_{|x-\Upsilon_2|\leq|x|/2}\Big]
&\leq\sum_{j=1}^N\mE \Big[ p_1 (x-\Upsilon_2)\1_{x-\Upsilon_2\in B_\eps(z_j)}\Big]\\
&\leq\sum_{j=1}^N\sup_{z\in B_\eps(z_j)}p_1 (z)\mP(|x-\Upsilon_2-z_j|\leq \eps)\\
&\leq  \sum_{j=1}^N\sup_{z\in B_\eps(z_j)}p_1 (z) \, \frac{c_1}{(1+|x-z_j|)^d \phi (1+|x-z_j|) }\\
&\leq  \frac{c_1c_4}{(1+|x|/2)^d \phi (1+|x|/2)} .
\end{align*}
 This together with  \eqref{e:1.1},  \eqref{ER6} and 
 Lemma \ref{Le22}(ii)  yields that 
 $$
p^{\kappa,\phi}_{0,1}(x) \leq \frac{c_5}{ (1+|x| )^d \phi (1+|x| )} \quad \hbox{for } x\in \mR^d.
 $$
 Combining with \eqref{e:2.17}, we get 
  the desired estimate  \eqref{ER11}.
 \end{proof}

\begin{remark}\label{R:2.7} \rm 
 The strong Markov process $X^{\kappa}_{t, s}$ of \eqref{e:2.3a} has infinitesimal
generator
$$
\sL^{\kappa}_t f(x) =  \int_{\R^d}
\left( f(x+z)-f(x)-\nabla f (x)\cdot z^{(\phi)} \right) \frac{\kappa (t, z)}{|z|^d \phi (|z|) } \dif z.
$$
Suppose
\begin{equation}\label{e:2.27}
\wt X^\kappa_{t,s}:=\int^s_t\!\!\!\int_{\mR^d}z \wt N(\dif r,\dif z)+\int^s_t\!\!\!\int_{\{z|>1\}} z\, \frac{\kappa(r,z)}{|z|^d \phi (|z|) }\dif z\dif r,
\end{equation} 
which has infinitesimal generator 
$$ 
\wt \sL^{\kappa}_t f(x) :=  \int_{\R^d}
\left( f(x+z)-f(x)-\nabla f (x)\cdot z \1_{\{|z|\leq 1\}} \right) \frac{\kappa (t, z)}{|z|^d \phi (|z|) } \dif z.
$$
Clearly,
\begin{equation}\label{e:2.28}
\wt X^\kappa_{t,s} =  X^\kappa_{t,s} + \int_t^s b(r) \dif r
\quad {\mbox{and}} \quad 
\wt \sL^{\kappa}_t =\sL^{\kappa}_t + b(r) \cdot \nabla, 
\end{equation} 
where  
\begin{equation}\label{e:2.29}
b(r)=
\begin{cases}
- \int_{\{|z|\leq 1\}} z \, \frac{\kappa (r, z)}{|z|^d \phi (|z|) } \dif z
&\hbox{ in Case$^\phi_1$ }, \\
0   &\hbox{ in Case$^\phi_2$  },  \\
\int_{\{|z|> 1\}} z \, \frac{\kappa (r, z)}{|z|^d \phi (|z|) } \dif z
&\hbox{ in Case$^\phi_3$}.
\end{cases}
\end{equation}
  Denote by $\wt p^{\kappa}_{t, s}(x)$ 
the density function of $\wt X^\kappa_{t,s}$.
Then by \eqref{e:2.28},
$$
 \wt p^{\kappa}_{t, s}(x) =  p^{\kappa}_{t, s} \left(x - \int_t^s b(r) \dif r \right).
$$
Thus under conditions \eqref{Ass1a},
one can get two-sided estimates on $\wt p^{\kappa}_{t, s}(x)$ from that of 
$p^{\kappa}_{t, s}(x) $. 
 \end{remark}

\subsection{Fractional derivative estimates of $p^{\kappa,\phi}_{0,1}$}
 
 In this subsection we show the fractional derivative estimates of $p^{\kappa,\phi}_{0,1}$, that will be used to construct the heat kernel with 
 variable coefficients by Levi's method.

\bl\label{Le24}
Let $\phi$ be as in \eqref{e:1.1} and $p(x):=p^{\kappa,\phi}_{0,1}(x)$. Define
 $$
\delta^{(1)}_p(x,z):=p(x+z)-p(x),\quad  \delta^{(2)}_p(x,z):=p(x+z)+p(x-z)-2p(x)
$$
and
$$
\delta^{(3)}_p(x,z):=p(x+z)-p(x)-z\cdot\nabla p(x).
$$
Under assumption \eqref{Ass1a}, there is a constant $c=c(\Theta_1 )>0$ 
such that for all $x,x',z\in\mR^d$ and 
$i=1,2,3$,
 \begin{align}
|\delta^{(i)}_p(x,z)|\leq c\,\Big(\1_{\{|z|>1\}}
 \big(\rho_\phi(x+z)+\1_{\{i=2\}}\rho_\phi(x-z)\big)+
(|z|^{i\wedge 2}\wedge  |z|^{(i-2)\vee 0} )
\rho_\phi(x)\Big),\label{ER222}
\end{align}
and
\begin{align}\label{ER202}
\begin{split}
& |\delta^{(i)}_p(x,z)-\delta^{(i)}_p(x',z)| \\
&\leq c \left(|x-x'|\wedge 1\right) 
\Big(\1_{\{|z|>1\}}\Big(\rho_\phi(x+ z)+\rho_\phi(x'+ z)  
  + \1_{\{i=2\}} \big( \rho_\phi(x- z)+\rho_\phi(x'- z) \big) \Big)  \\
& \hskip 1truein   + (|z|^{i\wedge 2}\wedge  |z|^{(i-2)\vee 0}  )
\big(\rho_\phi(x)+\rho_\phi(x')\big)\Big).
\end{split}
\end{align}
\el

\begin{proof}
We only present the proof of  \eqref{ER202} for  $i=3$. 
The proofs for the other two cases are similar. 
Notice that
\begin{align*}
|\delta^{(2)}_p(x,z)-\delta^{(2)}_p(x',z)|&\leq|z|^2\,|x-x'|\int_{[0,1]^3}|\nabla^3 p(x+\theta_1 z+\theta_1\theta_2 z+\theta_3 (x-x'))|\dif\theta_1\dif\theta_2\dif\theta_3\\
&\stackrel{\eqref{ER101}}{\lesssim}|z|^2\,|x-x'|\int_{[0,1]^3}\rho_\phi(x+\theta_1 z+\theta_1\theta_2 z+\theta_3(x-x'))\dif\theta_1\dif\theta_2\dif\theta_3.
\end{align*}
If $|x|>4$, $|x-x'|\leq 1$ and $|z|\leq 1$, then due to $|x+\theta_1 z+\theta_1\theta_2 z+\theta_3(x-x')|\asymp |x|,$ we have
$$
\rho_\phi(x+\theta_1 z+\theta_1\theta_2 z+\theta_3(x-x'))\lesssim \rho_\phi(x).
$$
If $|x|\leq 4$, $|x-x'|\leq 1$ and $|z|\leq 1$, then
$$
\rho_\phi(x+\theta_1 z+\theta_1\theta_2 z+\theta_3(x-x'))\leq 1\lesssim \rho_\phi(x).
$$
Hence, for $|x-x'|\leq 1$ and $|z|\leq 1$, 
$$
|\delta^{(2)}_p(x,z)-\delta^{(2)}_p(x',z)|\lesssim|z|^2\,|x-x'|\,\rho_\phi(x).
$$
Similarly, if $|x-x'|>1$ and $|z|\leq 1$, then 
$$
|\delta^{(2)}_p(x,z)-\delta^{(2)}_p(x',z)|\leq |\delta^{(2)}_p(x,z)|+|\delta^{(2)}_p(x',z)|\lesssim|z|^2\,\left(\rho_\phi(x)+\rho_\phi(x')\right);
$$
if $|x-x'|\leq 1$ and $|z|>1$, then 
$$
|\delta^{(2)}_p(x,z)-\delta^{(2)}_p(x',z)|\lesssim|x-x'|\left(\rho_\phi(x+z)+(|z|+1)\rho_\phi(x)\right);
$$
if $|x-x'|>1$ and $|z|>1$, then 
$$
|\delta^{(2)}_p(x,z)-\delta^{(2)}_p(x',z)|\lesssim\rho_\phi(x+z)+\rho_\phi(x'+z)+(|z|+1)(\rho_\phi(x')+\rho_\phi(x)).
$$
Combining the above cases, we obtain \eqref{ER202} for  $i=3$.
\end{proof}

Using Lemma \ref{Le24}, it is easy to derive the following from definition, which is a counterpart of  \cite[Theorem 2.4]{Ch-Zh}.

\bl\label{Th24}
Let $\phi$ be as in \eqref{e:1.1} and $z^{(\phi)}$ be defined by \eqref{e:1.3}.  
  Let   $\Delta^{(\phi)}_p$
be defined by any one of \eqref{e:Delta} 
(regardless whether $\kappa (t, z)$ is symmetric in $z$ or not), 
and $p(x):=p^{\kappa,\phi}_{0,1}(x)$. 
Under condition \eqref{Ass1a}, there is a constant $c=c(\Theta_1 )>0$ such that for all $x\in\mR^d,$
\begin{align}
\int_{\mR^d} 
|\Delta^{(\phi)}_p(x,z)|
\, \dif z\leq c\rho_\phi(x), 
\label{e:2.25}
\end{align}
and for all $x,x'\in\mR^d$,
\begin{align}\label{ER7}
\int_{\mR^d}  |\Delta^{(\phi)}_p(x,z)-\Delta^{(\phi)}_p(x',z)|  \, \dif z
\leq  c \left( |x-x'|\wedge 1 \right)  \left( \rho_{\phi}(x)+\rho_\phi(x')\right).
\end{align}
\el

\begin{proof}
 We first consider the case that 
   $\Delta^{(\phi)}_p$  
 is defined by the first expression of \eqref{e:Delta}.
 By Lemma \ref{Le24}, it is easy to see that
$$
 |\Delta^{(\phi)}_p(x,z)|\leq c
\, \frac{ \1_{\{|z|>1\}} \rho_\phi(x+z)+
 \gamma^{(1)}_\phi(|z|)
\rho_\phi(x)}{|z|^d \phi (z)},
$$
where 
 $\gamma^{(1)}_\phi(r)$ 
is defined by \eqref{e:1.21}, and
$$
|\Delta^{(\phi)}_p(x;z)-\Delta^{(\phi)}_p(x';z)|
 \leq c \left(|x-x'|\wedge 1\right)
 \frac{\1_{\{|z|>1\}}\big(\rho_\phi(x+ z)+\rho_\phi(x'+ z)\big)+
  \gamma^{(1)}_\phi(|z|)
 (\rho_\phi(x)+\rho_\phi(x'))} {|z|^d \phi (z)}. 
 $$
Thus, to prove \eqref{e:2.25} and \eqref{ER7}, it suffices to show
\begin{align*}
I_1:=\int_{\mR^d}\frac{
  \gamma^{(1)}_\phi(|z|)
}{|z|^d\phi(|z|)}\dif z<\infty \quad \hbox{and} \quad 
I_2(x):=\int_{|z|>1}\frac{\rho_\phi(x+z)}{|z|^d\phi(|z|)}\dif z\lesssim\rho_\phi(x).
\end{align*}
Clearly, by definition 
$$
I_1=c\int^\infty_0\frac{
  \gamma^{(1)}_\phi(r)
}{r\phi(r)}\dif r<\infty.
$$
For $|x|>2$, we have
\begin{align*}
I_2 (x)&\leq \rho_\phi(x/2)\int_{1<|z|\leq|x|/2}\frac{\dif z}{|z|^d\phi(|z|)}
+\int_{|z|>|x|/2}\frac{\rho_\phi(x+z)}{|z|^d\phi(|z|)}\dif z\\
&\lesssim \rho_\phi(x)\int_{|z|>1}\frac{\dif z}{|z|^d\phi(|z|)}
+\frac{1}{|x|^d\phi(|x|)}\int_{\mR^d}\rho_\phi(z)\dif z\lesssim\rho_\phi(x),
\end{align*}
 which together with $\sup_{|x|\leq 2} I_2(x)<\infty$ yields  $I_2 \lesssim\rho_\phi$ on $\mR^d$.

When 
 $\Delta^{(\phi)}_p$  
is defined by the second expression of \eqref{e:Delta} (here we do not need to assume
that $\kappa (t, z)$ is symmetric in $z$),  we can establish \eqref{e:2.25}
and \eqref{ER7} in a similar way as above.
This completes the proof of the lemma. 
\end{proof}

\subsection{Continuous dependence of $p^{\kappa,\phi}_{0,1}(x)$ with respect to $\kappa$}

In this subsection we show the continuous dependence of $p^{\kappa,\phi}_{0,1}(x)$ in the point-wise sense with respect to $\kappa$.
\bl\label{Le26}
Let $\kappa_1$ and $\kappa_2$ be two kernels satisfying \eqref{Ass1a} with the same constant $\kappa_0$. 
Let $\phi$ be as in \eqref{e:1.1} and 
$$
p_1(x):=p^{\kappa_1,\phi}_{0,1}(x), \quad  p_2(x):=p^{\kappa_2,\phi}_{0,1}(x).
$$ 
There exists a constant 
$c=c(\Theta_1 )>0$
such that for all $x\in\mR^d$,
\begin{align}\label{ER308}
|\nabla^j p_1(x)-\nabla^j p_2(x)|\leq c\|\kappa_1-\kappa_2\|_\infty\rho_\phi(x) \quad 
\hbox{for } j=0,1, 
\end{align}
\begin{align}\label{ER77}
\int_{\mR^d}|\Delta^{(\phi)}_{p_1}(x,z)-\Delta^{(\phi)}_{p_2}(x,z)|\dif z\leq c\|\kappa_1-\kappa_2\|_\infty \, \rho_\phi(x),
\end{align}
   where $\Delta^{(\phi)}_{p_i}(x,z)$ is defined by any one of \eqref{e:Delta} regardless 
   whether $\kappa (t, z)$ is symmetric in $z$ or not. 
\el

\begin{proof}
Noticing that by \eqref{Ch} and \eqref{Den1},
$$
q_\lambda(x):=
p^{\lambda\kappa_1+(1-\lambda)\kappa_2}_{0,1}(x) 
=\int_{\mR^d}\e^{-\mathrm{i}x\cdot\xi}
\exp\left(\int_{\mR^d}(\e^{\mathrm{i}\xi\cdot z}-1-\mathrm{i}\xi\cdot z^{(\phi)})
\frac{\lambda\kappa_1(z)+(1-\lambda)\kappa_2(z)}{|z|^d\phi(|z|)}\dif z\right) \dif\xi,
$$
where $\kappa_1(z):=\int^1_0\kappa_1(r,z)\dif r$ 
and $\kappa_2(z):=\int^1_0\kappa_2(r,z)\dif r$.
We claim that 
\begin{equation}\label{e:2.26}
\p_\lambda q_\lambda(x)=(\sL^{\kappa_1}-\sL^{\kappa_2})q_\lambda(x).
\end{equation} 
Indeed, since
$$
\widehat {\sL^{\kappa_i} f}(\xi)=\left(\int_{\mR^d}(\e^{\mathrm{i}\xi\cdot z}-1-\mathrm{i}\xi\cdot z^{(\phi)})
\frac{\kappa_i(z)}{|z|^d\phi(|z|)}\dif z\right)\hat f(\xi)
$$
and
$$
\hat q_\lambda(\xi)=\exp\left( \int_{\mR^d}(\e^{\mathrm{i}\xi\cdot z}-1-\mathrm{i}\xi\cdot z^{(\phi)})
\frac{\lambda\kappa_1(z)+(1-\lambda)\kappa_2(z)}{|z|^d\phi(|z|)}\dif z\right) ,
$$
we have
$$
\p_\lambda \hat q_\lambda(\xi)=\widehat {\sL^{\kappa_1} q_\lambda}(\xi)-\widehat{\sL^{\kappa_2}q_\lambda}(\xi).
$$
By the uniqueness of Fourier transform, we get \eqref{e:2.26}.

 By the definition of $q_\lambda(x)$ and \eqref{e:2.26}, 
\begin{align*}
|p_1(x)-p_2(x)|
&=\left|\int^1_0\p_\lambda q_\lambda(x)\dif \lambda\right|
=\left|\int^1_0\!\!\!\int_{\mR^d}
\Delta^{(\phi)}_{q_\lambda}(x,z)(\kappa_1(z)-\kappa_2(z))
\dif z\dif \lambda\right|\\
&\leq \|\kappa_1-\kappa_2\|_\infty\int^1_0\!\!\!\int_{\mR^d} 
\big|\Delta^{(\phi)}_{q_\lambda}(x,z)\big|\dif z.
\end{align*}
Thus we obtain by \eqref{e:2.25} that 
\begin{align}\label{UY0}
|p_1(x)-p_2(x)|\leq c\|\kappa_1-\kappa_2\|_\infty\rho_\phi(x),
\end{align}
which establishes \eqref{ER308} for $j=0$. 

We next  use the convolution technique to show \eqref{ER308} for $j=1$ and \eqref{ER77}.
Let 
$$
p_0(x):=p^{\kappa_0/2,\phi}_{0,1}(x),\quad 
\bar p_i(x) :=p^{\kappa_i-\kappa_0/2,\phi}_{0,1}(x)   \   \hbox{ for }  i=1,2.
$$ 
Then we have
$$
p_i(x)=\int_{\mR^d}p_0(x-y)
\bar p_i (y) \dif y     \   \hbox{ for }  i=1,2. 
$$
We have by \eqref{ER101} and \eqref{UY0} applied to $\bar p_1$ and $\bar p_2$ that 
\begin{align*}
|\nabla p_1(x)-\nabla p_2(x)|&\leq \int_{\mR^d}|\nabla p_0(x-y)|\cdot
| \bar p_1(y)- \bar p_2(y)|\dif y\\
&\lesssim \|\kappa_1-\kappa_2\|_\infty\int_{\mR^d}\rho_\phi(x-y)\cdot \rho_\phi(y)\dif y
\lesssim \|\kappa_1-\kappa_2\|_\infty\rho_\phi(x).
\end{align*}
Similarly, 
by \eqref{e:2.25} and \eqref{UY0},
\begin{align*} 
\int_{\mR^d}|\Delta^{(\phi)}_{p_1}(x,z)-\Delta^{(\phi)}_{p_2}(x,z)|\dif z  
&\leq \int_{\mR^d}\left(\int_{\mR^d}|\Delta^{(\phi)}_{p_0}(x-y,z)|\dif z\right) 
| \bar p_1(y)- \bar p_2(y)|
\dif y\\
&\lesssim\|\kappa_1-\kappa_2\|_\infty\int_{\mR^d}\rho_\phi(x-y)\cdot \rho_\phi(y)\dif y
\lesssim \|\kappa_1-\kappa_2\|_\infty\rho_\phi(x).
\end{align*}
This completes the proof of the  Lemma.
\end{proof}

\section{Basic convolution inequalities}\label{S:3}

We first  list
some important properties about slowly varying functions.

\bp\label{Pr32}
 \begin{enumerate}[\rm (i)]
 \item Let $\ell\in\sS_0$ with the convention that $\ell(t)=\ell(1)$ for $t\geq 1$. 
 For any $\delta>0$, there is a constant $C=C(\ell, \delta)\geq 1$ such that for all $s,t>0$,
\begin{align}
\label{e:3.1} 
\frac{\ell(s)}{\ell(t)}\leq C\max\left\{\left(\frac{s}{t}\right)^\delta,\left(\frac{s}{t}\right)^{-\delta}\right\}.
\end{align}
 
\item
If $\ell\in\sS_0\cap\sD_0$, then 
 $\Gamma_\ell (\cdot ) =\int^\cdot_0\ell(s)/s \, \dif s\in\sS_0$ and
\begin{align}\label{Kaa}
 \lim_{t\to 0}\Gamma_\ell(t)/\ell(t)=\infty.
\end{align}
If $\ell\in\sR_\alpha$ for some $\alpha>0$, then 
 \begin{align}\label{Ka}
 \lim_{t\to 0}\Gamma_\ell(t)/\ell(t)=1/\alpha.
\end{align}

\smallskip

\item Let $g:\mR_+\to\mR_+$ be a continuous function with $g(0)=0$. Suppose that
$$
0<\varliminf_{t\to 0}g(\lambda t)/g(t)\leq \varlimsup_{t\to 0}g(\lambda t)/g(t)<\infty,\quad \lambda>0.
$$
For any $\ell\in\sS_0$, we have $\ell\circ g\in\sS_0$.

\smallskip

\item For $\ell\in\sS_0\cap\sD_0$  
and an increasing positive function $\phi$ on $[0, \infty)$ that satisfies \eqref{e:1.1},
 we have $\ell_\phi\in\sS_0\cap\sD_0$.
 \end{enumerate}
\ep

\begin{proof}
(i) Estimate \eqref{e:3.1} 
follows by Potter's theorem (see \cite[(ii) of Theorem 1.5.6]{Bi-Go-Te}). 
\medskip\\
(ii) It follows by \cite[Proposition 1.5.9b and Proposition  1.5.10]{Bi-Go-Te}.
\medskip\\
(iii) For any $0<\lambda_0<\lambda_1<\infty$, by UCT theorem for slowly varying functions (see \cite[Theorem 1.2.1]{Bi-Go-Te}), it holds that
$$
\lim_{t\to 0}\sup_{\lambda\in[\lambda_0,\lambda_1]}|\ell(\lambda t)/\ell(t)-1|=0.
$$
Fix $\lambda>0$. By the assumption, there are $t_0>0$ and interval $[\lambda_0,\lambda_1]$ such that 
$$
g(\lambda t)/g(t)\in[\lambda_0,\lambda_1],\quad t\in(0,t_0].
$$
Hence,
$$
\lim_{t\to 0}\ell(g(\lambda t))/\ell(g(t))=1.
$$
(iv) By \eqref{e:1.1} we have
\begin{align}\label{SA1}
\left(\frac{R}{c_2^\phi r}\right)^{1/\beta_2}\leq\frac{\phi^{-1}  (R)}{\phi^{-1}  (r)}\leq \left(\frac{R}{c^\phi_1r}\right)^{1/\beta_1} \quad \hbox{for }  0<r\leq R\leq 1.
\end{align}
By (iii), we immediately have $\ell_\phi=\ell\circ\phi^{-1}\in\sS_0$. Moreover, the above estimate   implies
$$
\phi^{-1}(s)\leq (c_2^\phi s)^{1/\beta_2},
$$
 and so, by the increase of $\ell$ and a change of variable,
 \begin{align*}
\Gamma_{\ell_\phi}(t)&=\int^t_0\ell\circ\phi^{-1}(s)/s\dif s
\leq \int^t_0\ell\left((c_2^\phi s)^{1/\beta_2}\right)/s\dif s
=\beta_2\Gamma_\ell\Big((c_2^\phi t)^{1/\beta_2}\Big).
\end{align*}
The proof is complete.
 \end{proof}

The following lemma   
 plays a crucial role in  our construction of the heat kernel of $\sL^\kappa_t $   by Levi's method. 
 Recall the definition of $\ell_\phi$ from \eqref{e:1.14} 
 and  $\rho_\phi (t, x)$  from \eqref{Def8}.

\bl\label{Le33}
Let $\phi$ be as in \eqref{e:1.1} and $\rho_\phi(t,x)$ as in \eqref{Def8}.
 For a function $\ell$ on $\mR_+$, define
\begin{align}\label{HH}
h^\ell_\phi(t,x):=\ell\left(\phi^{-1} (t)+|x|\right)\rho_\phi(t,x).
\end{align}
\begin{enumerate}[\rm (i)]
 \item If $\ell\in\cup_{\alpha\in[0,1)}\sR_\alpha$, then there is  some $C=C(\ell,\Theta_1)>0$ so that   
\begin{align}\label{DK2}
\int_{\mR^d}h^\ell_\phi(t,x)\dif x\leq C\frac{\ell_\phi(t)}{t}
\quad \hbox{for all } t\in(0,1]. 
\end{align}

\item For any $\ell_1,\ell_2\in\cup_{\alpha\in[0,1)}\sR_\alpha$, 
there is a constant $C=C(\ell_1,\ell_2,\Theta_1 )\geq 1$ 
such that for all $0<s<t<\infty$ and $x,y\in\mR^d$,
\begin{align}\label{DK1}
\int_{\mR^d}h_\phi^{\ell_1}(t-s,x-y)h_\phi^{\ell_2}(s,y)\dif y\leq C\left(\frac{\ell_\phi(t-s)}{t-s}+\frac{\ell_\phi(s)}{s}\right)h_\phi^{\ell_1\vee\ell_2}(t,x).
\end{align}
\end{enumerate}
\el

\begin{proof}
(i) By the definition, \eqref{EK0} and the use of polar coordinate, we have
\begin{align}
\int_{\mR^d}h^\ell_\phi(t,x)\dif x&\lesssim\int^\infty_0\ell\left(\phi^{-1} (t)+r\right)\left(\phi^{-1}(t)+r\right)^{-d}\left(\phi \big(\phi^{-1}(t)+r\big)\right)^{-1}r^{d-1}\dif r \nonumber \\
&\leq \int^\infty_0 
\frac{\ell \big(\phi^{-1} (t)+r\big)}{  \big( \phi^{-1}(t)+r \big) \phi \big(\phi^{-1}(t)+r\big)} \dif r 
 =\int^\infty_{\phi^{-1} (t)}\frac{\ell(r)} {r\phi(r)}\dif r  \nonumber \\
& = \ell(1) \int^\infty_1\frac{1} {r\phi(r)}\dif r+\int^1_{\phi^{-1} (t)} \frac{\ell(r)} {r\phi(r)}\dif r. 
\label{e:3.7}
\end{align}
 Here we have used the convention that $\ell(r)=\ell(1)$ for $r\geq 1$.
 By the lower bound of \eqref{e:1.1} and 
\eqref{e:3.1} 
with $\delta = \beta_1/2$,  for $t\in (0, 1]$, 
\begin{align*}
\int^1_{\phi^{-1} (t)} \frac{\ell(r)} {r\phi(r)}\dif r
&= \frac{\ell (\phi^{-1}(t)) }{t}  \int^1_{\phi^{-1} (t)} \frac{\ell(r)} {\ell (\phi^{-1}(t)) } 
\frac{\phi (\phi^{-1}(t))} {r\phi(r)}\dif r \\
&\lesssim  \frac{\ell (\phi^{-1}(t)) }{t}  \int^1_{\phi^{-1} (t)} \frac{r^\delta} {\phi^{-1}(t)^\delta } 
\frac{ \phi^{-1}(t)^{\beta_1}} {r^{1+\beta_1}}\dif r   \lesssim   \frac{\ell (\phi^{-1}(t)) }{t} . 
\end{align*}
Taking $\delta = \beta_1/2$ again  in 
\eqref{e:3.1}, 
we have by \eqref{e:1.1} that  $t\lesssim \ell (\phi^{-1}(t))$
for $t\in (0, 1]$. The above display together with \eqref{e:3.7} and \eqref{e:1.1a} 
 yields the desired estimate \eqref{DK2}.

\medskip

(ii) Without loss of generality, we assume $\ell_1,\ell_2\in\sS_0$. By the definition of slowly varying function, it is easy to see that
$\ell_1\vee\ell_2\in\sS_0$. By 
\eqref{e:3.1}
with $\delta=1$, there is a constant $c_0>0$ such that for all $u<w$,
$$
\frac{u^d (\ell_1\vee\ell_2)(w)}{w^d(\ell_1\vee\ell_2)(u)}\leq c_0\left(\left(\frac{u}{w}\right)^{d+1}\vee\left(\frac{u}{w}\right)^{d-1}\right) \leq c_0.
$$
Thus, if we let $\gamma(u):=u^d\phi(u)$, then by the increase of $\phi$, we have 
for all $\lambda\geq 1$ and $0<u<\lambda w$,
\begin{align}\label{HP1}
 \frac{\gamma(u)}{(\ell_1\vee\ell_2)(u)}\leq c_0\frac{\gamma(\lambda w)}{(\ell_1\vee\ell_2)(\lambda w)}\leq c_\lambda\frac{\gamma(w)}{(\ell_1\vee\ell_2)(w)}.
\end{align}
On the other hand, as 
$$
u\vee w\leq u+w\leq 2(u\vee w)  \quad \hbox{for }   u,w>0,
$$
we have by 
\eqref{e:3.1} 
with $\delta=1$ again that  for all $u,w>0$,
$$
\frac{(\ell_1\vee\ell_2)(u\vee w)}{(\ell_1\vee\ell_2)(u+w)}\leq c_0\max\left\{\frac{u+w}{u\vee w},\frac{u\vee w}{u+w}\right\}\leq 2 c_0.
$$
This  together with $\gamma(u+w)\leq \gamma(2(u\vee w))\leq 2^{d+\beta_2}c^\phi_2 \gamma(u\vee w)$  yields
\begin{align}\label{HP2}
\frac{\gamma(u+w)}{(\ell_1\vee\ell_2)(u+w)}\lesssim\left(\frac{\gamma}{\ell_1\vee\ell_2}\right)(u\vee w)
\leq\frac{\gamma(u)}{\ell_1(u)}+\frac{\gamma(w)}{\ell_2(w)}.
\end{align}
Now let $g(t,x):=\phi^{-1} (t)+|x|$. Since for $0<s<t<\infty$ and $x,y\in\mR^d$,
$$
g(t,x)\lesssim g(t-s,x-y)+g(s,y),
$$
we have by \eqref{HP1} and \eqref{HP2},  
\begin{align*}
\frac{\gamma(g(t,x))}{(\ell_1\vee\ell_2)(g(t,x))}
&\lesssim \left(\frac{\gamma}{(\ell_1\vee\ell_2)}\right)(g(t-s,x-y)+g(s,y))\\
&\lesssim\frac{\gamma(g(t-s,x-y))}{\ell_1(g(t-s,x-y))}+\frac{\gamma(g(s,y))}{\ell_2(g(s,y))}.
\end{align*}
Hence,
\begin{align*}
&\frac{\ell_1(g(t-s,x-y))}{\gamma(g(t-s,x-y))}\times\frac{\ell_2(g(s,y))}{\gamma(g(s,y))}\\
&\quad\lesssim\left(\frac{\ell_1(g(t-s,x-y))}{\gamma(g(t-s,x-y))}+\frac{\ell_2(g(s,y))}{\gamma(g(s,y))}\right)
\times\frac{(\ell_1\vee\ell_2)(g(t,x))}{\gamma(g(t,x))},
\end{align*}
which together with \eqref{EK0} yields 
$$
h_\phi^{\ell_1}(t-s,x-y)h_\phi^{\ell_2}(s,y)\lesssim\left(h_\phi^{\ell_1}(t-s,x-y)+h_\phi^{\ell_2}(s,y)\right)h_\phi^{\ell_1\vee\ell_2}(t,x).
$$
Integrating both sides in $y$ and by \eqref{DK2}, we obtain \eqref{DK1}.
\end{proof}

\section{Proof of Theorem \ref{heat}}\label{S:4}

Now we consider the space and time dependent nonlocal operator
$\sL^{\kappa}_t$ defined by \eqref{Non}, 
 and give a proof for Theorem \ref{heat}.  
  For each fixed $y\in\mR^d$, let $\sL^{\kappa_y}_t$ be the freezing operator
\begin{align}\label{Def7} 
\sL^{\kappa_y}_t f(x)=\int_{\mR^d} \Delta^{(\phi)}_f(x,z)\,\kappa(t,y,z)\dif z,
\end{align}
where $\Delta^{(\phi)}_f$ is the difference operator defined by \eqref{e:Delta}. 
 Let $p^{(y)}_{t,s}(x):=p^{\kappa_y}_{t,s}(-x)$ 
be the heat kernel of operator $\sL^{\kappa_y}_t$
as given by \eqref{Den1}. 
Equivalently,  $p^{(y)}_{t,s}(x)$ is the probability transition density of the time-inhomogenous L\'evy process
associated with $\sL^{\kappa_y}_t$ starting from position $x$ at time $t$ to be at the origin 0 at time $s$.
We know from  \eqref{DL6} that it satisfies  
\begin{align}
\p_t p^{(y)}_{t,s}(x)+\sL^{\kappa_y}_t p^{(y)}_{t,s}(x)=0 \  \hbox{for } s>t
\hbox{ with } \lim_{t\uparrow s}p^{(y)}_{t,s}(x)=
\delta_{\{0\}} (x),  \label{ES2}
\end{align}
 where  $\delta_{\{ 0\}} (x)$ denotes the usual Dirac 
measure concentrated at the origin $0$. 
 
Following Levi's idea, we   seek   heat kernel $p^\kappa_{t,s} (x,y)$ of $\sL^{\kappa}_t$ of the following form:
\begin{align}\label{ER65}
p^\kappa_{t,s} (x,y)= p^{(y)}_{t,s}(x-y)+\int^s_t\!\!\!\int_{\mR^d}p^{(z)}_{t,r}(x-z)q_{r,s}(z,y)\dif z\dif r,
\end{align}
where $q_{r,s}(z,y)$  is  some suitable function   to be determined. 
If the above $p^\kappa_{t,s} (x,y)$ is a heat kernel for $ \sL^{\kappa}_t$,
that is, for each $t<s$ and $x,y\in\mR^d$,
$$
\p_t p^\kappa_{t,s}(x ,y) +\sL^{\kappa}_t p^\kappa_{t,s}(\cdot ,y) (x) =0,
$$
formally differentiate both sides of \eqref{ER65} with respect to $t$ would  yield
$$
\sL^{\kappa}_tp^\kappa_{t,s} (\cdot ,y) (x)= 
 \big(\sL^{\kappa_y}_tp^{(y)} _{t,s}\big) (x -y) 
+ q_{t,s}(x, y) +\int^s_t\!\!\!\int_{\mR^d} 
\big(\sL^{\kappa_z}_t p^{(z)}_{t,r}\big) (x-z)
q_{r,s}(z,y)\dif z\dif r. 
$$
Applying $\sL^{\kappa}_t$ on both sides of \eqref{ER65} in $x$-variable formally gives
$$
\sL^{\kappa}_tp^\kappa_{t,s} (\cdot ,y) (x) = 
\big(\sL^{\kappa_x}_t p^{(y)}_{t,s}\big) (x -y) 
+\int^s_t\!\!\!\int_{\mR^d}  
\big(\sL^{\kappa_x}_t p^{(z)}_{t,r}\big) (x-z) 
q_{r,s}(z,y)\dif z\dif r. 
$$
Subtracting the above two displays and defining 
\begin{align}\label{WM3}
q^{(0)}_{t,s}(x,y):=   \left(\sL^{\kappa_x}_{t}-\sL^{\kappa_y}_{t} \right)p^{(y)}_{t,s} (\cdot)  (x-y),
\end{align}
we conclude that $q_{t,s}(x,y)$ must satisfy 
\begin{align}\label{EU2}
q_{t,s}(x,y)=q^{(0)}_{t,s}(x,y)+\int^s_t\!\!\!\int_{\mR^d}q^{(0)}_{t,r}(x,z)q_{r,s}(z,y)\dif z\dif r
\end{align}
for any $t<s$ and $x, y \in \mR^d$.
For $n\in\mN$, define $q^{(n)}_{t,s}(x,y)$ recursively by
\begin{align}\label{EU22} 
q^{(n)}_{t,s}(x,y):=\int^s_t\!\!\!\int_{\mR^d}q^{(0)}_{t,r}(x,z)q^{(n-1)}_{r,s}(z,y)\dif z\dif r. 
\end{align}
Iterating the identity \eqref{EU2}  repeatedly, we get for $N\geq 1$,
\begin{equation}\label{e:4.7}
q_{t,s}(x,y) =\sum_{n=0}^N  q^{(n)}_{t,s}(x,y)+\int^s_t\!\!\!\int_{\mR^d}q^{(N)}_{t,r}(x,z)q_{r,s}(z,y)\dif z\dif r.
\end{equation} 
If the remainder tends to zero as $N\to \infty$, we would get 
 \begin{equation}\label{e:4.8} 
 q_{t,s}(x,y):=\sum_{n=0}^\infty  q^{(n)}_{t,s}(x,y).
 \end{equation} 
Our approach is that, instead of showing the remainder in \eqref{e:4.7} tends to zero, we   show that the  infinite sum 
in the above  
converges absolutely and locally uniformly   and  the function  $q_{t,s}(x,y)$ defined 
by \eqref{e:4.8} satisfies the integral equation \eqref{e:4.7},   using the estimates  
obtained in the last two sections.
We then establish  rigorously that 
 the function $p^\kappa_{t, s}(x, y)$ defined by \eqref{ER65} in terms of 
 $q_{t,s}(x,y)$ of \eqref{e:4.8} is indeed a time-inhomogeneous heat kernel for $\sL^{\kappa}_t$ 
 with desired regularity and estimates.  The positivity of $q_{t,s}(x,y)$  and the uniqueness of the heat kernel
 are obtained through a maximum principle for non-local operator $\sL^\kappa_t$.  
  
   {\it
    Throughout the remaining of this section, we assume $\ell \in \sS_0\cap \sD_0$
 and  one of the following holds:
 \begin{enumerate} 
\item[{\bf (H1)}] If $\kappa(t,x,z)=\kappa(t,x,-z)$, we assume \eqref{e:1.1}, 
 ({\bf A$^{(0)}_\phi$})
and \eqref{Con1}.
\item[{\bf (H2$'$)}] If $\kappa(t,x,z)\not=\kappa(t,x,-z)$, we assume \eqref{e:1.1}, 
 ({\bf A$^{(1)}_\phi$}),
\eqref{Con1}
 and \eqref{Con2}.
  \end{enumerate}
 }
 
\subsection{Solving integral equation \eqref{EU2}}
Our first step is to 
 show that the function $q_{t,s}(x, y)$ given by \eqref{e:4.8} solves the integral equation \eqref{EU2}. 
We will use scaling to reduce the consideration of heat kernel $p^{\kappa}_{t, s}$
to the case of $t=0$ and $s=1$. 
 For this, we define for each fixed  $t<s$ with $s-t\leq 1$, 
\begin{align}\label{DL1}
\begin{split}
\wt\kappa_y(r,z):=\kappa(t+r(s-t), y, \phi^{-1}(s-t)z),\quad 
\wt \phi (u)&:= \phi (u\phi^{-1}(s-t))/ (s-t).
\end{split}
\end{align}
Note that $\wt \phi $ satisfies \eqref{e:1.1} with the same constants $c_1^\phi$, $c_2^\phi$ and $0<\beta_1\leq \beta_2<\infty$, 
 and, in view of 
   \eqref{Con0},
\begin{equation}\label{e:4.10} 
   \int^\infty_0\frac{\gamma^{(i)}_{\wt\phi}(r)}{r \wt\phi(r)}\dif r\leq \cA^{(i)}_\phi<\infty
 \quad \hbox{for }  i=0,1.
\end{equation} 
  Clearly, $\wt\kappa_y(r,z)$ satisfies \eqref{Con2}, and  we have  by \eqref{Con1} that
\begin{equation}\label{Con1'} 
\kappa^{-1}_0\leq \wt\kappa_y (r, z )\leq \kappa_0, \quad
 | \wt\kappa_x (r, z)- \wt\kappa_y (r,z)|
  \leq\ell^2\left(|x-y|\right).
\end{equation}
By Proposition \ref{Pr21}, we have for any $x,y$ and $t<s$,
\begin{align}\label{EK5}
q^{(0)}_{t,s}(x,y)=
\frac{\left(\Big(\sL^{\wt\kappa_x,\wt\phi}_0-\sL^{\wt\kappa_y,\wt\phi}_0\Big)p^{\wt\kappa_y, \wt \phi}_{0,1} (\cdot) \right)\Big((x-y)/\phi^{-1}(s-t)\Big)}
{(s-t) \left( \phi^{-1} (s-t)\right)^{d }}.
\end{align}
Noticing that by definition \eqref{Def7}, \eqref{e:2.25} and \eqref{Con1'},
\begin{align*}
\left|\Big(\sL^{\wt\kappa_x,\wt\phi}_0-\sL^{\wt\kappa_y,\wt\phi}_0\Big)p^{\wt\kappa_y, \wt \phi}_{0,1}\right|(z)
 \lesssim\ell^2(|x-y|)\,\rho_{\wt\phi}(z),
\end{align*}
and also by definition \eqref{Def8} and \eqref{DL1},
\begin{align}\label{LS1}
\rho_{\wt\phi}((x-y)/\phi^{-1}(s-t))=(s-t) \left( \phi^{-1} (s-t)\right)^d\rho_\phi(s-t,x-y),
\end{align}
we conclude that there is a positive constant $C_0=C_0(\Theta_1, \cA^{(0)}_\phi)$ so that for  any $(t,x;s,y)\in\mD^1_0$,
\begin{align}\label{EK8}
\begin{split}
|q^{(0)}_{t,s}(x,y)|&\leq C_0\ell^2(|x-y|)\,\rho_{\phi}(s-t,x-y)\leq C_0h^{\ell^2}_\phi(s-t,x-y),
\end{split}
\end{align}
where $h^{\ell^2}_\phi$ is defined by \eqref{HH}.

\medskip

The following theorem extends  \cite[Theorem 3.1]{Ch-Zh} to  the time-dependent  and mixed 
stable-like non-local 
operator setting
of this paper, and relaxes the H\"older continuous assumption on $x\mapsto \kappa (t, x, z)$ to Dini continuity. 
Recall the definition of $\ell_\phi$ from \eqref{e:1.14}.

\bt\label{T:4.1}
 Let $q^{(n)}_{t,s}(x,y)$ be defined as  by \eqref{WM3} and \eqref{EU22}. 
 Under either  {\bf (H1)} or  {\bf (H2$'$)}, 
 there is an $\eps_0>0$ such that
the series $q_{t,s}(x,y):=\sum_{n=0}^{\infty}q^{(n)}_{t,s}(x,y)$ is
 absolutely and locally uniformly convergent on $\mD^{\eps_0}_0$
and solves the integral equation  \eqref{EU2}. Moreover, for each $t>0$,
 $(s,x,y)\mapsto q_{t,s}(x,y)$ is jointly continuous in $\mK^{\eps_0}_t:=(t,t+\eps_0)\times\mR^d\times\mR^d$, 
 and has the following estimates:
there is a constant $c_1=c_1 (\eps_0, \Theta)>0$ so that on $\mD^{\eps_0}_0$,
\begin{align}
|q_{t,s}(x,y)| \leq c_1h^{\ell^2}_\phi(s-t,x-y)\leq c_1\|\ell\|_\infty h^\ell_\phi(s-t,x-y), \label{eq3}
\end{align}
where $h^\ell_\phi$ is defined by \eqref{HH}, and
\begin{align}\label{eq4}
|q_{t,s}(x,y)-q_{t,s}(x',y)| \leq c_1\frac{\ell(|x-x'|)}{\ell_\phi(s-t)} \left(h^\ell_\phi(s-t,x-y)+h^\ell_\phi(s-t,x'-y)\right).
\end{align}
\et

\begin{proof}
(i) Let $C_2:=2C_1C_0$, where $C_1$ is the constant in \eqref{DK1} 
 associated to $\ell_1=\ell_2=\ell^2$,
 and $C_0$ is from \eqref{EK8}.
We use induction method to show that for all $(t,x;s,y)\in\mD^1_0$,
$$
|q^{(n)}_{t,s}(x,y)|\leq C^{n+1}_2\left(\Gamma_{\ell^2_\phi}(s-t)\right)^nh^{\ell^2}_\phi(s-t,x-y),
$$
where $\Gamma_{\ell^2_\phi}(t):=\int^t_0\ell^2_\phi(s)\dif s/s$.
First of all, for $n=0$, it is true by \eqref{EK8}.
Suppose now that it has been proven for some $n\in\mN$. Then 
by \eqref{EK8}, the induction hypothesis and  \eqref{DK1}, we have
\begin{align*}
|q^{(n+1)}_{t,s}(x,y)|
&\leq C_2^{n+1}\left(\Gamma_{\ell^2_\phi}(s-t)\right)^nC_0
\int^s_t\!\!\int_{\mR^d}h^{\ell^2}_\phi(r-t,x-z)h^{\ell^2}_\phi(s-r,z-y)\dif z\dif r\\
&\leq C^{n+1}_2\left(\Gamma_{\ell^2_\phi}(s-t)\right)^nC_0\left(C_1
\int^s_t\left(\frac{\ell^2_\phi(r-t)}{r-t}+\frac{\ell^2_\phi(s-r)}{s-r}\right) \dif r\right)h^{\ell^2}_\phi(s-t,x-y)\\
&=C^{n+1}_2\left(\Gamma_{\ell^2_\phi}(s-t)\right)^nC_0 \left(2C_1\Gamma_{\ell^2_\phi}(s-t)\right)h^{\ell^2}_\phi(s-t,x-y)\\
&= C_2^{n+2}\left(\Gamma_{\ell^2_\phi}(s-t)\right)^{n+1} h^{\ell^2}_\phi(s-t,x-y).
\end{align*}
  Now by (iv) of Proposition \ref{Pr32},
we can choose $\eps_0\in(0,1)$ small enough so that
$$
 \Gamma_{\ell^2_\phi}(\eps_0)\leq 1/(2C_2).
$$
Thus $q_{t,s}(x,y)=\sum_{n=0}^{\infty}q^{(n)}_{t,s}(x,y)$ converges absolutely and locally uniformly on $\mD^{\eps_0}_0$ and
$$
|q_{t,s}(x,y)|\leq\sum_{n=0}^{\infty}|q^{(n)}_{t,s}(x,y)|\leq 2C_2h^{\ell^2}_\phi(s-t,x-y).
$$
(ii) By \eqref{e:2.3a}, \eqref{DH1} and \eqref{Den1}, one sees that $\{(t,s,x)\mapsto p^{(y)}_{t,s}(x): y\in\mR^d\}$ is equi-continuous in any compact subsets of 
$\{(t,s,x)\in\mR_+\times\mR_+\times\mR^d: t<s\}$. On the other hand, by \eqref{ER308} and Proposition \ref{Pr21}, it is easy to see that
$y\mapsto p^{(y)}_{t,s}(x)$ is continuous for each $t<s$ and $x\in\mR^d$. Hence, $(t,x;s,y)\mapsto p^{(y)}_{t,s}(x-y)$ is continuous on $\mD^1_0$.
Moreover, by the definition of $q^{(0)}$, one sees that for each $t>0$, $(s,x,y)\mapsto q^{(0)}_{t,s}(x,y)$ 
is continuous on $(t,t+1)\times\mR^d\times\mR^d$.
Furthermore, by definition \eqref{EU22} and induction method, for each $n\in\mN$,
$(s,x,y)\mapsto q^{(n)}_{t,s}(x,y)$ is continuous on $(t,t+1)\times\mR^d\times\mR^d$. So $(s,x,y)\mapsto q_{t,s}(x,y)$ is continuous on 
$\mK^{\eps_0}_t$.
\\
\\
(iii) In this step we show that for all $(t,x;s,y),(t,x';s,y)\in\mD^1_0$,
\begin{align}\label{WM1}
|q^{(0)}_{t,s}(x,y)-q^{(0)}_{t,s}(x',y)|\lesssim \frac{\ell(|x-x'|)}{\ell_\phi(s-t)} \left(h^{\ell^2}_\phi(s-t,x-y)+h^{\ell^2}_\phi(s-t,x'-y)\right).
\end{align}
First of all, if $|x-x'|\geq\phi^{-1}(s-t)$, then by \eqref{EK8} and the increase of $\ell$, 
we clearly have the above estimate.
Next we assume
$$
|x-x'|\leq\phi^{-1}(s-t).
$$
Write $p(x):=p^{\wt\kappa_y, \wt \phi}_{0,1}(x)$. 
 Recall the definition of $\Delta^{(\wt\phi)}_p$ from \eqref{e:Delta}. 
We have by definition \eqref{DL1}, Lemma \ref{Th24} and \eqref{Con1'},
 \begin{align*}
&\left|\Big(\sL^{\wt\kappa_x,\wt\phi}_0-\sL^{\wt\kappa_y,\wt\phi}_0\Big)p(z_1)
-\Big(\sL^{\wt\kappa_{x'},\wt\phi}_0-\sL^{\wt\kappa_y,\wt\phi}_0\Big)p(z_2)\right|\\
&\leq\int_{\mR^d}  
\left|\Delta^{(\wt\phi)}_{p}(z_1,z)-\Delta^{(\wt\phi)}_{p}(z_2,z)\right|
\, |\wt\kappa_x(0,z)-\wt\kappa_y(0,z)|\, \dif z
+\int_{\mR^d} 
 \left|\Delta^{(\wt\phi)}_{p}(z_2,z)\right|
\, |\wt\kappa_x(0,z)-\wt\kappa_{x'}(0,z)|\, \dif z\\
&\lesssim \ell^2(|x-y|)(|z_1-z_2|\wedge 1)\left(\rho_{\wt\phi}(z_1)+\rho_{\wt\phi}(z_2)\right)+\ell^2(|x-x'|)\rho_{\wt\phi}(z_2),
\end{align*}
 where the implicit constant depends on  $\Theta $, 
   especially on $\cA^{(0)}_\phi$ (resp. $\cA^{(1)}_\phi$) in the symmetric (resp. non-symmetric) case
   of $z\mapsto \kappa (t, x, z)$. 
 By \eqref{EK5} and \eqref{LS1} and taking $z_1=\frac{x-y}{\phi^{-1}(s-t)}$ and $z_2=\frac{x'-y}{\phi^{-1}(s-t)}$, we get for $|x-x'|\leq\phi^{-1}(s-t)$,
\begin{align}\label{HH8}
\begin{split}
&|q^{(0)}_{t,s}(x,y)-q^{(0)}_{t,s}(x',y)|\lesssim
\ell^2(|x-y|)\frac{|x-x'|}{\phi^{-1}(s-t)}\rho_{\phi}(s-t,x-y)\\
&\qquad+\left(\ell^2(|x-y|)\frac{|x-x'|}{\phi^{-1}(s-t)}+\ell^2\left(|x-x'|\right)\right)\rho_{\phi}(s-t,x'-y).
\end{split}
\end{align}
Since $\ell\in\sS_0\cap\sD_0$ is bounded,  by Proposition \ref{Pr32} we have
$$
\tfrac{r}{R}\lesssim \tfrac{\ell(r)}{\ell(R)}\lesssim 1\mbox{ for }r\leq R\ \mbox{ and }\ \ell(s+t)\lesssim\ell(s)+\ell(t).
$$
Thus, by the definition \eqref{HH} of $h^{\ell^2}_\phi$, from \eqref{HH8}, we immediately have \eqref{WM1}.
\medskip\\
(iv) By \eqref{WM1}, \eqref{eq3} and Lemma \ref{Le33}, we have
\begin{align}\label{WM2}
\begin{split}
&\int^s_t\!\!\int_{\mR^d}|q^{(0)}_{t,r}(x,z)-q^{(0)}_{t,r}(x',z)|\,|q_{r,s}(z,y)|\dif z\dif r\\
&\lesssim \ell(|x-x'|)\int^s_t\frac{1}{\ell_\phi(r-t)}\left[\frac{\ell^2_\phi(r-t)}{r-t}+\frac{\ell^2_\phi(s-r)}{s-r}\right]\dif r \\
&\quad\times\left(h^{\ell^2}_\phi(s-t,x-y)+h^{\ell^2}_\phi(s-t,x'-y)\right).
\end{split}
\end{align}
Clearly, we have
$$
\int^s_t\frac{\ell_\phi(r-t)}{r-t}\dif r=\int^{s-t}_0\frac{\ell_\phi(r)}{r}\dif r=\Gamma_{\ell_\phi}(s-t)<\infty.
$$
Write
$$
\int^s_t\frac{\ell^2_\phi(s-r)}{\ell_\phi(r-t) (s-r)}\dif r=\left(\int^s_{(s+t)/2}+\int^{(s+t)/2}_t\right)\frac{\ell^2_\phi(s-r)}{\ell_\phi(r-t) (s-r)}\dif r=:I_1+I_2.
$$
For $I_1$, since $\ell_\phi\in\sS_0\cap\sD_0$, we have
$$
I_1\lesssim\frac{1}{\ell_\phi((s-t)/2)}\int^{(s-t)/2}_0\frac{\ell^2_\phi(r)}{r}\dif r\leq\int^{(s-t)/2}_0\frac{\ell_\phi(r)}{r}\dif r=\Gamma_{\ell_\phi}((s-t)/2).
$$
For $I_2$, since $s\mapsto s/\ell_\phi(s)\in\sR_1$ by $1/\ell_\phi\in\sS_0$, we have by \eqref{Ka} that
$$
I_2\lesssim\frac{\ell^2_\phi(s-t)}{s-t}\int^{(s-t)/2}_0\frac{\dif r}{\ell_\phi(r)}\lesssim\ell_\phi(s-t).
$$
Combining these with \eqref{WM2}, \eqref{WM1} and \eqref{EU2}, 
 we obtain \eqref{eq4}.
\end{proof}

\br \rm 
In order to obtain estimate \eqref{eq4}, we need to borrow some regularity from the spatial variable to compensate 
the time singularity (see \eqref{WM2}). This is the only reason that we have to assume \eqref{Con1} for the square of some Dini's function. 
\er

\bc\label{Co42} 
Suppose   either assumption  {\bf (H1)} or  {\bf (H2$'$)} holds. Let 
 $p^\kappa_{t,s}(x,y)$ be defined by \eqref{ER65} and $\eps_0$ be as in Theorem \ref{T:4.1}. Then $p^\kappa_{t,s}(x,y)$ is continuous 
on $\mD^{\eps_0}_0$ and there are constants $c_0,c_1>0$ and $\delta>0$ such that 
\begin{align}\label{DL7}
p^\kappa_{t,s}(x,y)\leq c_0\,(s-t)\rho_\phi(s-t,x-y)  \quad \mbox{on } \mD^{\eps_0}_0,
\end{align}
and if $|x-y|\leq\phi^{-1}(s-t)\leq\delta$, then
\begin{align}\label{DL8}
p^\kappa_{t,s}(x,y)\geq c_1\phi^{-1}(s-t)^{-d}.
\end{align}
\ec
\begin{proof}
First of all, by Proposition \ref{Pr21} and Theorem \ref{Th1}, there is a constant $c_2>1$ such that on $\mD^1_0$,
\begin{align}\label{DL99}
c^{-1}_2(s-t)\rho_\phi(s-t,x-y)\leq p^{(y)}_{t,s}(x,y)\leq c_2(s-t)\rho_\phi(s-t,x-y).
\end{align}
Define
\begin{align}\label{DH2}
\tilde\ell(t):=\int^{\phi(t)}_0\ell\circ\phi^{-1}(s)/s\dif s=\Gamma_{\ell_\phi}\circ\phi(t).
\end{align}
By Proposition \ref{Pr32}, we know that $\ell_\phi\in\sS_0\cap\sD_0$ and $\tilde\ell\in\sS_0$. Moreover, by \eqref{Kaa},
\begin{align}\label{DH3}
\ell_\phi(t)\lesssim\Gamma_{\ell_\phi}(t)=\tilde\ell\circ\phi^{-1}(t)=\tilde\ell_\phi(t),\quad t\in[0,1].
\end{align}
Thus, by \eqref{eq3}, \eqref{DL99} and \eqref{DK1}, we have
\begin{align*}
&\int^s_t\!\!\!\int_{\mR^d}|p^{(z)}_{t,r}(x,z)q_{r,s}(z,y)|\dif z\dif r\lesssim
\int^s_t\!\!\!\int_{\mR^d}(r-t)\rho_\phi(r-t,x-z)h^\ell_\phi(s-r,z-y)\dif z\dif r\\
&\qquad\qquad\leq\int^s_t\frac{r-t}{\tilde\ell_\phi(r-t)}\left(\int_{\mR^d}h_\phi^{\tilde\ell}(r-t,x-z)h^\ell_\phi(s-r,z-y)\dif z\right)\dif r\\
&\qquad\qquad\lesssim h_\phi^{\tilde\ell}(s-t,x-y)\int^s_t\frac{r-t}{\tilde\ell_\phi(r-t)}\left(\frac{\tilde\ell_\phi(r-t)}{r-t}+\frac{\ell_\phi(s-r)}{s-r}\right)\dif r\\
&\qquad\qquad\lesssim h_\phi^{\tilde\ell}(s-t,x-y) \left((s-t)+\int^s_t\frac{r-t}{\tilde\ell_\phi(r-t)}\frac{\ell_\phi(s-r)}{s-r}\dif r\right).
\end{align*}
Note that by \eqref{DH3} and \eqref{DH2},
\begin{align*}
&\int^s_t\frac{r-t}{\tilde\ell_\phi(r-t)}\frac{\ell_\phi(s-r)}{s-r}\dif r
=\left(\int^{(s+t)/2}_t+\int^s_{(s+t)/2}\right)\left(\frac{r-t}{\tilde\ell_\phi(r-t)}\frac{\ell_\phi(s-r)}{s-r}\right)\dif r\\
&\qquad\lesssim\frac{\ell_\phi(s-t)}{s-t}\int^{s-t}_0\frac{r}{\tilde\ell_\phi(r)}\dif r+\frac{s-t}{\tilde\ell_\phi(s-t)}\int^{s-t}_0\frac{\ell_\phi(r)}{r}\dif r\\
&\qquad\lesssim\frac{\tilde\ell_\phi(s-t)}{s-t}\int^{s-t}_0\frac{r}{\tilde\ell_\phi(r)}\dif r+(s-t)\lesssim s-t,
\end{align*}
where the last step is due to $s\mapsto s^2/\tilde\ell_\phi(s)\in\sR_2$ and \eqref{Ka}.
Hence,
\begin{align}\label{DL5}
\int^s_t\!\!\!\int_{\mR^d}|p^{(z)}_{t,r}(x,z)q_{r,s}(z,y)|\dif z\dif r\leq c_3(s-t)h_\phi^{\tilde\ell}(s-t,x-y),
\end{align}
which together with \eqref{ER65} and \eqref{DL9} yields \eqref{DL7}.

On the other hand, if $|x-y|\leq\phi^{-1}(s-t)$, then by \eqref{DL99} and \eqref{DL5},
$$
p^\kappa_{t,s}(x,y)\geq (s-t)\rho_\phi(s-t,x-y)\left(c^{-1}_2-c_3\tilde\ell_\phi(s-t)\right).
$$
Choosing $\delta$ be small enough, we get \eqref{DL8}.

Finally, since $(t,x;s,y)\mapsto p^{(y)}_{t,s}(x-y)$ is continuous on $\mD^1_0$, by \eqref{ER65}, Theorem \ref{T:4.1} 
and the dominated convergence theorem, 
one sees that $(t,x;s,y)\mapsto p^\kappa_{t,s}(x,y)$ is continuous on $\mD^{\eps_0}_0$
\end{proof}

\subsection{Gradient and fractional derivative estimates of $p^\kappa_{t,s}$}
This section is similar to \cite[Sections 3.2 and 3.3]{Ch-Zh}. We only point out the main points.
The following lemma follows easily from \eqref{DH1}, \eqref{Den1}, Theorem \ref{Th1} and equation \eqref{ES2}.

\bl\label{Le41}
Suppose  either   {\bf (H1)} or  {\bf (H2$'$)} holds. 
For  each $j\in\mN$, $s>0$ and $y\in\mR^d$, the mapping $(t,x)\mapsto\nabla^j p^{(y)}_{t,s} (x-y)$ is continuous on $[0,s)\times\mR^d$.
Moreover, there is a constant $C>0$ such that for all $(t,x;s,y)\in\mD^1_0$,
\begin{align}
\label{e:4.26}
|\nabla^j p^{(y)}_{t,s} (x-y)|\leq \frac{ C(s-t)}{ (\phi^{-1}(s-t) )^j}  \, \rho_\phi(s-t,x-y)
\end{align}
and
\begin{align}
\lim_{t\uparrow s}\sup_{x\in\mR^d}\left|\int_{\mR^d}p^{(y)}_{t,s}(x-y)\dif y-1\right|=0.\label{DL10}
\end{align}
\el
\begin{proof}
The estimate \eqref{e:4.26} follows from \eqref{Sca}, \eqref{ER101} and \eqref{LS1}.
We next show \eqref{DL10}. 
By \eqref{ER308}, \eqref{Sca}, \eqref{Con1'}, \eqref{LS1} and \eqref{DK2}, we have
\begin{align*}
\left|\int_{\mR^d}p^{(y)}_{t,s}(x-y)\dif y-1\right|&=\left|\int_{\mR^d}\left(p^{(y)}_{t,s}(x-y)-p^{(x)}_{t,s}(x-y)\right)\dif y\right|\\
 &\lesssim (s-t)\int_{\mR^d}\ell^2(|x-y|) \rho_\phi(s-t,x-y)\dif y\\
&\leq(s-t)\int_{\mR^d}h^{\ell^2}_\phi(s-t,x-y)\dif y\lesssim \ell^2_\phi(s-t),
\end{align*}
 where the implicit constant $C$ is independent of $x$ and $s-t$.
Thus we get \eqref{DL10}.
\end{proof}

To show the gradient and fractional derivative estimates, by \eqref{ER65} we write 
\begin{align}\label{FB1}
\begin{split}
p^\kappa_{t,s} (x,y) = & \,  p^{(y)}_{t,s}(x-y)
+\int^s_{\frac{s+t}{2}}\!\int_{\mR^d} p^{(z)}_{t,r}(x-z)q_{r,s}(z,y)\dif z\dif r\\
&+\int^{\frac{s+t}{2}}_t\!\!\!\int_{\mR^d} p^{(z)}_{t,r}(x-z)q_{r,s}(z,y)\dif z\dif r \\
 =: & \, \sum_{i=1}^3J_i(t,x;s,y).
\end{split}
\end{align}
Recall from \eqref{e:1.14} that 
$\ell_\phi (t):=\ell (\phi^{-1} (t))$ and $\Gamma_{\ell}(t):=\int^t_0\frac{\ell(s)}{s}\dif s$.

\bl\label{Le46}
Assume that in addition to the assumption  {\bf (H1)} or  {\bf (H2$'$)},  condition \eqref{Con3} holds as well. 
Then for each $0\leq t<s$ with $s-t\leq 1$ and $y\in\mR^d$, 
$ p^\kappa_{t,s} (x,y)$ is continuously differentiable in $x\in \mR^d$.  Moreover, 
there is a constant $C>0$ such that for all $(t,x;s,y)\in\mD^1_0$,
\begin{align}\label{DL11}
|\nabla_x p^\kappa_{t,s} (x,y)|\leq C\left(\frac{s-t}{\phi^{-1}(s-t)}+M^\phi_\ell\circ\phi^{-1}(s-t)\right)\rho_\phi(s-t,x-y),
\end{align}
where $M^\phi_\ell (t)$ is the function defined by \eqref{Con3}.
\el

\begin{proof}
For $J_1(t,x;s,y)$ in \eqref{FB1}, we have by \eqref{e:4.26}
$$
|\nabla J_1(t,\cdot;s,y)(x)|\lesssim  \frac{s-t}{\phi^{-1}(s-t)} \, \rho_\phi(s-t,x-y).
$$
For $J_2(t,x;s,y)$ in \eqref{FB1}, we have by \eqref{DL9}, \eqref{eq3} and \eqref{DK1}, 
\begin{align*}
|\nabla J_2(t,\cdot;s,y)(x)|&\lesssim \int^s_{\frac{s+t}{2}}\!\int_{\mR^d}\frac{r-t}{\phi^{-1}(r-t)}\rho_\phi(r-t,x-z)h^{\ell^2}_\phi(s-r,z-y)\dif z\dif r\\
&\lesssim \rho_\phi(s-t,x-y)\int^s_{\frac{s+t}{2}}\frac{r-t}{\phi^{-1}(r-t)}\left(\frac{1}{r-t}+\frac{\ell_\phi(s-r)}{s-r}\right)\dif r\\
&\lesssim \frac{s-t}{\phi^{-1}(s-t)}\Big(1+\Gamma_{\ell_\phi}(s-t)\Big)\rho_\phi(s-t,x-y).
\end{align*}
For $J_3(t,x;s,y)$, we approximate it by 
\begin{align}\label{BV1}
J^{(\eps)}_3(t,x;s,y):=\int^{\frac{s+t}{2}}_{t+\eps}\!\!\!\int_{\mR^d} p^{(z)}_{t,r}(x-z)q_{r,s}(z,y)\dif z\dif r,\ \ \eps\in(0,\tfrac{s-t}{2}).
\end{align}
Fix $\eps\in(0,\tfrac{s-t}{2})$. By \eqref{e:4.26} and \eqref{eq3}, we can exchange $\nabla$ with the integral and arrive at
\begin{align}
\nabla J^{(\eps)}_3(t,\cdot;s,y)(x)&=\int^{\frac{s+t}{2}}_{t+\eps}\!\!\!\int_{\mR^d} \nabla p^{(z)}_{t,r}(x-z)q_{r,s}(z,y)\dif z\dif r\no\\
&=\int^{\frac{s+t}{2}}_{t+\eps}\!\!\!\int_{\mR^d}\nabla  p^{(z)}_{t,r}(x-z)(q_{r,s}(z,y)-q_{r,s}(x,y))\dif z\dif r\no\\
&\quad+\int^{\frac{s+t}{2}}_{t+\eps}\left(\int_{\mR^d}\left(\nabla p^{(z)}_{t,r}-\nabla p^{(x)}_{t,r}\right)(x-z)\dif z\right)q_{r,s}(x,y)\dif r\no\\
&=:K^{(\eps)}_1(t,x;s,y)+K^{(\eps)}_2(t,x;s,y),\no 
\end{align}
where in the second equality we have used
$$
\int_{\mR^d}\nabla p^{(x)}_{t,r}(x-z)\dif z=0,
$$
For $K^{(\eps)}_1(t,x;s,y)$,  by \eqref{eq4} and \eqref{DK2}, \eqref{DK1}, we have
\begin{align*}
&\quad|K^{(\eps)}_1(t,x;s,y)|\\
&\lesssim
\int^{\frac{s+t}{2}}_{t+\eps}\!\!\!\int_{\mR^d}\frac{r-t}{\phi^{-1}(r-t)}\rho_\phi(r-t,x-z)\frac{\ell(|x-z|)}{\ell_\phi(s-r)}
\Big(h^\ell_\phi(s-r,z-y)+h^\ell_\phi(s-r,x-y)\Big)\dif z\dif r\\
&\lesssim\frac{1}{\ell_\phi(s-t)}
\int^{\frac{s+t}{2}}_t\!\!\!\frac{r-t}{\phi^{-1}(r-t)}\!\int_{\mR^d}h^\ell_\phi(r-t,x-z)\Big(h^\ell_\phi(s-r,z-y)+h^\ell_\phi(s-r,x-y)\Big)\dif z\dif r\\
&\lesssim\left(\frac{h^\ell_\phi(s-t,x-y)}{\ell_\phi(s-t)}\right)
\int^{\frac{s+t}{2}}_t\left(\frac{r-t}{\phi^{-1}(r-t)}\right)\left(\frac{\ell_\phi(r-t)}{r-t}+\frac{\ell_\phi(s-r)}{s-r}\right)\dif r\\
 &\lesssim\left(\frac{1}{\ell_\phi(s-t)}\int^{s-t}_0\frac{\ell_\phi(r)}{\phi^{-1}(r)}\dif r+\frac{1}{s-t}\int^{s-t}_0\frac{r}{\phi^{-1}(r)}\dif r\right)\rho_\phi(s-t,x-y)\\
&=\left(\int^{\phi^{-1}(s-t)}_0\left( \frac{\ell(r)}{\ell_\phi(s-t) r}+\frac{\phi(r)}{(s-t)r}\right) 
     \dif\phi(r)\right)\rho_\phi(s-t,x-y)\\
&=M^\phi_\ell\circ\phi^{-1}(s-t) \, \rho_\phi(s-t,x-y).
\end{align*}
For $K^{(\eps)}_2(t,x;s,y)$, noting that by \eqref{ER308},
$$
\big|\nabla p^{\tilde\kappa_y, \wt \phi}_{0,1}-\nabla p^{\tilde\kappa_x, \wt \phi}_{0,1}\big|(z)
\lesssim\|\tilde\kappa_x-\tilde\kappa_y\|_\infty \, \rho_{\wt\phi}(z)\lesssim\ell(|x-y|) \, \rho_{\wt\phi}(z),
$$
by Proposition \ref{Pr21} and   \eqref{LS1}, we have
\begin{align*}
|K^{(\eps)}_2(t,x;s,y)|&\lesssim\int^{\frac{s+t}{2}}_t \frac{r-t}{\phi^{-1}(r-t)}\left(\int_{\mR^d}h^\ell_\phi(r-t,x-z)\dif z\right)h^\ell_\phi(s-r,x-y)\dif r\\
&\lesssim\left(\int^{\frac{s+t}{2}}_t \frac{\ell_\phi(r-t)}{\phi^{-1}(r-t)}\dif r\right)h^\ell_\phi(s-t,x-y)
\lesssim\left(\int^{s-t}_0 \frac{\ell_\phi(r)}{\phi^{-1}(r)}\dif r\right)\rho_\phi(s-t,x-y)\\
&=\left(\int^{\phi^{-1}(s-t)}_0 \frac{\ell(r)}{r}\dif \phi(r)\right)\rho_\phi(s-t,x-y)\lesssim M^\phi_\ell\circ\phi^{-1}(s-t)\rho_\phi(s-t,x-y),
\end{align*}
where the above implicit constant is independent of $\eps$.
Moreover, from the above proof, it is also easy to see that
$$
\lim_{\eps\downarrow 0}K^{(\eps)}_i(t,x;s,y)=K^{(0)}_i(t,x;s,y),\ \ i=1,2,
$$
locally uniformly. Moreover, by the dominated convergence theorem, 
$$
x\mapsto K^{(0)}_i(t,x;s,y)\mbox{ is continuous.}
$$
As $J^{(\eps)}_3(t,x;s,y)$ converges to $J_3(t,x;s,y)$ pointwise and $\nabla_x J^{(\eps)}_3(t,x;s,y)=
 K^{(\eps)}_1(t,x;s,y) +  K^{(\eps)}_i(t,x;s,y)$, we conclude that 
 $J_3(t,x;s,y)$   is  differentiable in $ x $   and  $ \nabla_x J_3(t,x;s,y)= K^{(0)}_1(t,x;s,y) +  K^{(0)}_i(t,x;s,y)$,
 which is  continuous  in $x$. 
Summing the above up, we have shown that $ p^\kappa_{t,s} (x,y)$ is continuously differentiable in $x\in \mR^d$
whose gradient has the desired estimate \eqref{DL11}. 
\end{proof}

\bl\label{Le42}
Suppose either  {\bf (H1)} or  {\bf (H2$'$)} holds. 
For each $0\leq t<s$ and $y\in\mR^d$, the mapping $x\mapsto\sL^{\kappa_x}_t p^{(y)}_{t,s} (x)$ is continuous
and for fixed $t_0<s$ and $x\in\mR^d$,
$$
\lim_{t\downarrow t_0}\left|\sL^{\kappa_x}_t p^{(y)}_{t,s} (x)-\sL^{\kappa_x}_t p^{(y)}_{t_0,s} (x)\right|=0.
$$
Moreover, there is a constant $C>0$ such that for all $(t,x;s,y)\in\mD^1_0$,
\begin{align}\label{DL9}
\int_{\mR^d} 
\big|\Delta^{(\phi)}_{p^{(y)}_{t,s}}(x-y,z)\big| \, \dif z 
\leq C\rho_\phi(s-t,x-y),
\end{align} 
where $\Delta^{(\phi)}_p$ is defined in \eqref{e:Delta}.
\el

\begin{proof}
  We first show \eqref{DL9}.  
By Proposition \ref{Pr21},    \eqref{e:2.25}  and \eqref{LS1}, 
we have
\begin{align*}
 \int_{\mR^d}\!\big|\Delta^{(\phi)}_{p^{(y)}_{t,s}}(x-y,z)\big| \, \dif z 
& = (s-t)^{-1} \left( \phi^{-1} (s-t)\right)^{-d }
\!\!\!\int_{\mR^d} 
 \left|\Delta^{(\phi)}_{p^{\tilde\kappa, \wt \phi}_{0,1}}\left((x-y)/\phi^{-1}(s-t),z\right)\right| \, \dif z \\
&  \lesssim  (s-t)^{-1} \left( \phi^{-1} (s-t)\right)^{-d }\rho_{\wt\phi}((x-y)/\phi^{-1}(s-t))
=\rho_\phi(s-t,x-y).
\end{align*}
This proves  \eqref{DL9}. 
 Note that by Remark \ref{R:1.1}(i), 
\begin{align*}
\sL^{\kappa_x}_t p^{(y)}_{t,s} (x)=
 \int_{\mR^d} \Delta^{(\phi)}_{p^{(y)}_{t,s}} (x, z)   \, \dif z.    
\end{align*}
The desired continuity follows by the dominated convergence theorem and Lemma \ref{Le41}.
\end{proof}

By Lemma \ref{Le42}, the following lemma can be proved in a similar way as that for Lemma \ref{Le46}.

\bl\label{Le45}
 Suppose the condition of Theorem \ref{heat} holds. 
 For each $0\leq t<s$ with $s-t\leq 1$ and $y\in\mR^d$, 
 $x\mapsto p^\kappa_{t,s} (x,y)$ is ponitwisely $\sL^\kappa_t$-differentiable in the sense that the integral in
\eqref{Non} and \eqref{Non1} is absolutely convergent for every $x\in \mR^d$. 
Moreover, $x\mapsto\sL^\kappa_t p^\kappa_{t,s} (\cdot ,y)(x)$
is continuous
and for fixed $t_0<s$ and $x\in\mR^d$,
$$
\lim_{t\downarrow t_0}\left|\sL^{\kappa}_t p^{\kappa}_{t,s} (x,y)-\sL^{\kappa}_t p^{\kappa}_{t_0,s} (x,y)\right|=0.
$$ 
 Furthermore, there is a constant $C>0$ such that for all $(t,x;s,y)\in\mD^1_0$,
\begin{align}\label{SA0}
 \int_{\mR^d}\big|\Delta^{(\phi)}_{p^\kappa_{t,s}}(x-y,z)\big| \, \dif z 
\leq C  \left(\frac{\Gamma_{\ell_\phi}(s-t)}{\ell_\phi(s-t)}\right)
\rho_\phi(s-t,x-y).
\end{align}
\el

\begin{proof}
 Recall from \eqref{FB1}, $ p^\kappa_{t,s} (x,y) = \sum_{k=1}^3  J_k(t,x;s,y)$. By \eqref{DL9},
$$
\int_{\mR^d}
\big|\Delta^{(\phi)}_{J_1(t,\cdot; s,y)}(x,z)\big|\, \dif z 
 \lesssim \rho_\phi(s-t,x-y).
$$
For $J_2(t,x;s,y)$, we have by \eqref{DL9}, \eqref{eq3} and \eqref{DK1}, 
\begin{align*}
\int_{\mR^d}\big|\Delta^{(\phi)}_{J_2(t,\cdot; s,y)}(x,z)\big|\, \dif z 
&\leq\int^s_{\frac{s+t}{2}}\!\int_{\mR^d}\int_{\mR^d}\big|\Delta^{(\phi)}_{p^{(\bar z)}_{t,s}}(x-\bar z,z)\big|\, \dif z
\,  |q_{r,s}(\bar z,y)|\dif \bar z\dif r\\
&\lesssim \int^s_{\frac{s+t}{2}}\!\int_{\mR^d}\rho_\phi(r-t,x-\bar z)h^\ell_\phi(s-r,\bar z-y)\dif \bar z\dif r\\
&\lesssim \rho_\phi(s-t,x-y)\int^s_{\frac{s+t}{2}}\left(\frac{1}{r-t}+\frac{\ell_\phi(s-r)}{s-r}\right)\dif r\\
&\lesssim \rho_\phi(s-t,x-y)\Big(1+\Gamma_{\ell_\phi}\big(\tfrac{s-t}2\big)\Big)
\lesssim \rho_\phi(s-t,x-y). 
\end{align*}
For $J_3(t,x;s,y)$ in   \eqref{FB1}, we can use the same approximation as used in Lemma \ref{Le46}. 
Here we omit approximation procedure  and only make the following calculations.
Note that
\begin{align}
\Delta^{(\phi)}_{J_3(t,\cdot; s,y)}(x,z)&=\int^{\frac{s+t}{2}}_{t}\!\!\!\int_{\mR^d} \Delta^{(\phi)}_{p^{(\bar z)}_{t,r}}(x-\bar z,z)q_{r,s}(\bar z,y)\dif \bar z\dif r\no\\
&=\int^{\frac{s+t}{2}}_{t}\!\!\!\int_{\mR^d} \Delta^{(\phi)}_{p^{(\bar z)}_{t,r}}(x-\bar z,z)(q_{r,s}(\bar z,y)-q_{r,s}(x,y))\dif \bar z\dif r\no\\
&\quad+\int^{\frac{s+t}{2}}_{t}\!\left(\int_{\mR^d} \left(\Delta^{(\phi)}_{p^{(\bar z)}_{t,r}}(x-\bar z,z)-\Delta^{(\phi)}_{p^{(x)}_{t,r}}(x-\bar z,z)\right)\dif\bar z\right)q_{r,s}(x,y)\dif r,
\end{align}
where we have used
$$
\int_{\mR^d} \Delta^{(\phi)}_{p^{(x)}_{t,r}}(x-\bar z,z)\dif\bar z=0.
$$
We therefore have
\begin{align}
\int_{\mR^d}\big|\Delta^{(\phi)}_{J_3(t,\cdot; s,y)}(x,z)\big| \, \dif z 
&\leq\int^{\frac{s+t}{2}}_t\!\!\!\int_{\mR^d}\left(\int_{\mR^d} \left|\Delta^{(\phi)}_{p^{(\bar z)}_{t,r}}(x-\bar z,z)\right|\, 
\dif z \right)|q_{r,s}(\bar z,y)-q_{r,s}(x,y)|\dif \bar z\dif r\no\\
&\quad+\int^{\frac{s+t}{2}}_t\!\!\!\int_{\mR^d}\dif\bar z\left(\int_{\mR^d} \left|\Delta^{(\phi)}_{p^{(\bar z)}_{t,r}-p^{(x)}_{t,r}}(x-\bar z,z)\right| \, \dif z \right)
|q_{r,s}(x,y)|\dif r\no\\
&=:K_1(t,x;s,y)+K_2(t,x;s,y).\label{BV2}
\end{align}
For $K_1(t,x;s,y)$, we have by \eqref{DL9} and \eqref{eq4},
\begin{align*}
|K_1(t,x;s,y)|&\lesssim
\int^{\frac{s+t}{2}}_t\!\!\!\int_{\mR^d}\rho_\phi(r-t,x-z)\frac{\ell(|x-z|)}{\ell_\phi(s-r)}\Big(h^\ell_\phi(s-r,z-y)+h^\ell_\phi(s-r,x-y)\Big)\dif z\dif r\\
&\lesssim\frac{1}{\ell_\phi(s-t)}
\int^{\frac{s+t}{2}}_t\!\!\!\int_{\mR^d}h_\phi^\ell(r-t,x-z)\Big(h^\ell_\phi(s-r,z-y)+h^\ell_\phi(s-r,x-y)\Big)\dif z\dif r\\
&\stackrel{\eqref{DK1}}{\lesssim}\left(\frac{h^\ell_\phi(s-t,x-y)}{\ell_\phi(s-t)}\right)
\int^{\frac{s+t}{2}}_t\left(\frac{\ell_\phi(r-t)}{r-t}+\frac{\ell_\phi(s-r)}{s-r}\right)\dif r\\
&\leq\frac{2\Gamma_{\ell_\phi}(s-t)}{\ell_\phi(s-t)}h^\ell_\phi(s-t,x-y)\lesssim\frac{\Gamma_{\ell_\phi}(s-t)}{\ell_\phi(s-t)}\rho_\phi(s-t,x-y),
\end{align*}
where the last inequality 
is due to the fact that $\ell (t)=\ell (1)$ for $t\geq 1$   so $\ell \lesssim 1$ on $[0, \infty)$ 
and consequently  $h^\ell_\phi(s-t,x-y)\lesssim  \rho_\phi(s-t,x-y)$. 
For $K_2(t,x;s,y)$, noting that by   Proposition \ref{Pr21} and \eqref{ER77},
\begin{align*}
\int_{\mR^d} \left|\Delta^{(\phi)}_{p^{(\bar z)}_{t,r}-p^{(x)}_{t,r}}(x-\bar z,z)\right|\, \dif z 
\lesssim  \rho_\phi (r-t, x-\bar z) \,   \ell(|x-\bar z|) \leq h^\ell_\phi (r-t, x-\bar z),
\end{align*} 
we have by     \eqref{eq3}, Lemma \ref{Le33}(i) and Proposition \ref{Pr32}(ii) that 
\begin{align*}
|K_2(t,x;s,y)| 
& \lesssim   \int^{\frac{s+t}{2}}_t  \left( \int_{\mR^d}h^\ell_\phi (r-t, x-\bar z) \dif \bar z\right)   h^\ell_\phi (s-r, x-y)    \dif r \\
& \lesssim   \left(  \int^{\frac{s+t}{2}}_t  \frac{\ell_\phi(r-t)}{r-t} \dif r \right) h^\ell_\phi (s-t, x-y)  \\
&\lesssim  \Gamma_\ell ((s-t)/2) \, \rho_\phi(s-t,x-y) \\
&\lesssim  \rho_\phi(s-t,x-y). 
\end{align*} 
 Combining the above calculations, we obtain \eqref{SA0}.
As for the desired continuity of $\sL^{\kappa}_t p^{\kappa}_{t,s} (x,y)$ in $t$, it follows from Theorem \ref{T:4.1},  \eqref{EU2}, Lemma \ref{Le42}, \eqref{eq3}, \eqref{DK1}  and the dominated convergence theorem. 
The proof is now complete.
\end{proof}

\subsection{A maximum principle} In this subsection we establish a maximum principle for operator $\sL^\kappa$,
which will be used to obtain  the uniqueness and positivity of heat kernels.

\bt\label{Le51} 
  For $T>0$, let $u(t,x)\in C_b([0,T)\times\mR^d)$ satisfy the following equation: for all $x\in\mR^d$ and Lebesgue almost all $t\in[0,T)$,
  $$
    \p_t u(t,x)+\sL^\kappa_tu(t,x)\leq 0,
    \quad \varliminf_{t\uparrow T}u(t,x)\geq 0.
  $$
  Assume that for each $t\in[0,T)$ and $x\in\mR^d$, 
  \begin{align}\label{AL2}
  \lim_{s\downarrow t}|\sL^\kappa_su(s,x)-\sL^\kappa_su(t,x)|=0,
  \end{align}
and that  in Case$^\phi_2$ and Case$^\phi_3$, when $\kappa (t, x, z)$ is not symmetric in $z$,  
  for each $t\in[0,T)$,
\begin{align}\label{AL1}
\mbox{$x\mapsto\nabla u(t,x)$ is continuous on $\mR^d$.}
\end{align}
Then we have
  \begin{align}\label{NM4}
    u(t,x)\geq 0,\quad (t,x)\in[0,T)\times\mR^d.
  \end{align}
\et

\begin{proof}
We only give the proof for the case when $\kappa (t, x, z)$ is not symmetric in $z$. 
The proof when  $\kappa (t, x, z)$ is   symmetric in $z$ is similar except that we use 
\eqref{Non1} for  the expression of $\sL^\kappa_t$.
  First of all, we assume that for all $t\in[0,T)$,
  $$
  \lim_{|x|\to\infty} u(t,x)=\infty
  $$
  and there is some constant $\delta <0$ so that 
    for each $x\in\mR^d$ and Lebesgue almost all $t\in[0,T)$,
  \begin{align}\label{HG1}
    \p_t u(t,x)+\sL^\kappa_tu(t,x)\leq \delta<0
  \end{align}
  Suppose that \eqref{NM4} is not true.
  Since $\lim_{|x|\to\infty} u(t,x)=\infty$ and $\varliminf_{t\uparrow T}u(t,x)\geq 0$,
  there must be a point $(t_0,x_0)\in [0,T)\times\mR^d$ such that
  $$
    u(t_0,x_0)=\inf_{(t,x)\in[0,T)\times\mR^d}u(t,x)<0.
  $$
In Case$^\phi_2$ and Case$^\phi_3$, since $x_0$ is 
a minimun 
point of $x\mapsto u(t_0,x)$, by \eqref{AL1} we have
  $$
    \nabla u(t_0,x_0)=0.
  $$
  Therefore, for each $s>0$,
  $$
    \sL^\kappa_s u(t_0,x_0)=\int_{\mR^d} \left( u(t_0,x_0+z)-u(t_0,x_0)-z^{(\phi)}\cdot\nabla u(t_0,x_0) \right)
    \frac{\kappa(s,x_0,z)}{|z|^d\phi(|z|)}\dif z\geq 0,
  $$
  and integrating both sides of \eqref{HG1}  from $t_0$ to $t$, we have
  \begin{align}
     u(t,x_0)-u(t_0,x_0) & \leq   (t-t_0)\delta-\int^t_{t_0}\sL^\kappa_s u(s,x_0)\dif s\no\\
    & \leq(t-t_0)\delta-\int^t_{t_0}\left(\sL^\kappa_su(s,x_0)-\sL^\kappa_su(t_0,x_0)\right)\dif s.\label{UT4}
  \end{align}
  Dividing both sides  by $t-t_0$ and letting $t\downarrow t_0$,  we obtain  by \eqref{AL2} 
  \begin{align*}
    0&\leq \delta+\lim_{t\downarrow t_0}\frac{1}{t-t_0}\int^t_{t_0}|\sL^\kappa_su(s,x_0)-\sL^\kappa_su(t_0,x_0)|\dif s=\delta<0,
    \end{align*}
  which is impossible. In other words, the infimum is achieved at the terminal time $T$, and \eqref{NM4} holds.

  Next, we  drop the restriction \eqref{HG1}.  For this, let
  $$
    f(x):=(1+|x|^2)^\alpha,\quad \alpha\in(0,\beta_1/2).
  $$
  For $\eps,\delta>0$, define
  $$
    u_{\eps,\delta}(t,x):=u(t,x)+\delta(T-t)+\eps\e^{-t} f(x).
  $$
  By easy calculations, one sees that for some $C>0$,
  $$
    |\sL^\kappa_t f(x)|\leq C(1+|x|^\alpha),
  $$
  and
  \begin{align*}
    \p_tu_{\delta,\eps}(t,x)+\sL^\kappa_t u_{\delta,\eps}(t,x)\leq-\delta+\eps\e^{-t}(\sL^\kappa_t f(x)-f(x))\leq-\delta/2<0,
  \end{align*}
  provided $\eps$ being small enough so that $\eps\e^{-t}(\sL_tf(x)-f(x))<\delta/2$. Clearly, 
  $$
    \lim_{x\to\infty}|u_{\eps,\delta}(t,x)|=\infty.
  $$
  Hence, by what we have proved,
  $$
    u_{\eps,\delta}(t,x)\geq 0.
  $$
  By letting $\eps\to 0$ and then $\delta\to 0$, we obtain \eqref{NM4}.
\end{proof}

\smallskip

\subsection{Proof of Theorem \ref{heat}}\label{S:4.4} 

Assume the conditions of Theorem \ref{heat} hold.
Let $\eps_0>0$ be defined as in Theorem \ref{T:4.1} and $p^\kappa_{t,s}(x,y)$ be defined by \eqref{ER65}.
For $f\in C^2_b(\mR^d)$, define 
$$
u(t,x):=P^\kappa_{t,s}f(x):=\int_{\mR^d}p^\kappa_{t,s}(x,y)f(y)\dif y.
$$

(i) It follows from  Lemma \ref{Le45} that  $x\to u(t, x)$  is 
 pointwisely
$\sL^\kappa_t$-differentiable in the sense that 
 the integrals in \eqref{Non} and \eqref{Non1} are absolutely convergent 
for every $x\in \mR^d$, 
and that for fixed $t_0<s$ and $x\in\mR^d$,
\begin{align}\label{EK3}
\lim_{t\downarrow t_0}\left|\sL^\kappa_t u(t,x)-\sL^\kappa_tu(t_0,x)\right|=0.
\end{align}

(ii) It follows from  Lemma \ref{Le46} that when $\kappa (t, x, z)$ is not symmetric, in Case$^\phi_2$ and Case$^\phi_3$
under condition \eqref{Con3},   
\begin{align}\label{EK2}
\mbox{$x\mapsto\nabla u(t,x)$ is continuous on $\mR^d$.}
\end{align}

(iii) For any bounded and uniformly continuous function $f$, by \eqref{DL10} and \eqref{DL5}, it is not hard to see that 
$$
\lim_{t\uparrow s}\|P^\kappa_{t,s}f-f\|_\infty=0.
$$
Moreover, by Lemma \ref{Le45},  Lemma \ref{Le46} and the discussion at the beginning of this section, one has 
that (see \cite{Ch-Zh} for more details)
\begin{align}\label{EK00}
u(t,x)=f(x)+\int^s_t\sL^\kappa_r u(r,x)\dif r,\  \ \forall (t,x)\in[0,s)\times\mR^d.
\end{align}
The maximum principle from Theorem \ref{Le51} gives 
the uniqueness of $p^\kappa_{t,s}(x,y)$ as well as the properties that 
\begin{align}\label{DL12}
p^\kappa_{t,s}(x,y)\geq 0\mbox{  and }\int_{\mR^d}p^\kappa_{t,s}(x,y)\dif y=1\ \mbox{ on $\mD^{\eps_0}_0$,}
\end{align}
and for all $0\leq t<r<s<\infty$ with $s-t\leq\eps_0$ and $x,y\in\mR^d$,
\begin{align}\label{eq211}
\int_{\mR^d}p^\kappa_{t,r} (x,z)p^\kappa_{r,s} (z,y)\dif z=p^\kappa_{t,s} (x,y).
\end{align}

Now we are in a position to give
\begin{proof}[{Proof of Theorem \ref{heat}}]
 Let $\eps_0>0$ be the constant from Theorem \ref{T:4.1}. 
 We have established in the above the existence and uniqueness of heat kernel $p^\kappa_{t,s}(x,y)$ on $\mD^{\eps_0}_0$ 
 that satisfies (i)-(iii) of Theorem \ref{heat} on $\mD^{\eps_0}_0$. 
 We now  extend the definition of $p^\kappa_{t,s}(x,y)$ and its properties in (i)-(iii)
   from $\mD^{\eps_0}_0$ to $\mD^\infty_0$ by \eqref{eq211} 
  as follows: If $\eps_0<s-t\leq 2\eps_0$, we define
  \begin{align}\label{Ext}
    p^\kappa_{t,s}(x,y)=\int_{\mR^d}p^\kappa_{t,\frac{t+s}{2}}(x,z)p^\kappa_{\frac{t+s}{2},s}(z,y)\dif z.
  \end{align}
  Proceeding this procedure, we can extend $p^\kappa$ to $\mD^\infty_0$ and the Chapman-Kolmogorov
   equation \eqref{eq211} holds for all
  $0\leq t<r<s<\infty$ and $x,y\in\mR^d$. In particular, equation \eqref{EQ} and (i), (ii), (iii) hold for  $p^\kappa_{t,s}(x,y)$.
  
  \smallskip
   Next  we show that   the heat kernel $p^\kappa_{t,s}(x,y)$ 
  enjoys properties {\bf (a)}-{\bf (f)}.

\smallskip
 \begin{enumerate} [\bf (a)]

\item  The upper bound estimate follows by \eqref{Ext}, \eqref{DL7} and \eqref{DK1}. Moreover, by \eqref{DL12} and \eqref{Ext} we also have
$$
p^\kappa_{t,s}(x,y)\geq 0.
$$
The lower bound will be proved in the next subsection.

\item  It follows by  \eqref{Ext}, \eqref{DL9} and \eqref{DK1}.

\item  It follows by  \eqref{Ext}, \eqref{DL11} and \eqref{DK1}.

\item  It follows by  \eqref{Ext} and \eqref{DL12}.

\item  It follows by  \eqref{Ext} and \eqref{eq211}.

\item  Fix $s>0$. Define 
$$
\tilde u(t,x):=f(x)+\int^s_t\!P^\kappa_{t,r}\sL^\kappa_rf(x)\dif r.
$$
By Fubini's theorem, it is easy to see that $\tilde u$ also satisfies equation \eqref{EK00} and \eqref{EK3}, \eqref{EK2}
(see \cite{Ch-Zh} and \cite{Ch-Hu-Xi-Zh}).
Thus by the maximum principle, we have $\tilde u(t,x)=u(t,x)=P^\kappa_{t,s}f(x)$. 
\end{enumerate}
This completes the proof of the theorem except for the lower bound in {\bf (a)} on the heat kernel  $p^\kappa_{t,s}(x,y)$ ,
which will be given separately in next subsection.
\end{proof}

\subsection{Proof of lower bound in \eqref{GR1}}

We know from the last subsection that 
$\{p^\kappa_{t,s}(x,y): (t,x;s,y)\in\mD^\infty_0\}$ is  
a family of transition probability density functions.
It uniquely determines a Feller process
$$ 
X:=\Big\{ \Omega, \, \sF, \,  (X_s)_{s\geq 0};  \, \mP_{t,x},  \, {(t,x)\in\mR_+\times\mR^d}\Big\}
$$
on $\mR^d$  with the property that
$$
  \mP_{t,x}\big(X_s=x,\,0\leq s\leq t\big)=1,
$$
and for $r\in[t,s]$ and $A\in\cB(\mR^d)$,
\begin{align}\label{Tr}
 \mE_{t,x}\left[ X_{s}\in A \,|\, \sF_{r}\right]
 =\int_{A}p^\kappa_{r,s}(X_{r},y)\dif y, 
\end{align}
where $\sF_s:=\sigma\{X_t, t\leq s\}$, $s\geq 0$, is the filtration generated by the Feller process $X$.
Moreover, for any $f\in C^2_b(\mR^d)$, it follows from (\ref{eqge}) and the Markov property of $X$ that under $\mP_{t,x}$, with respect to the filtration
$\{\sF_s; s\geq 0\}$
\begin{align}
  M^f_s:=f(X_s)-f(X_t)-\int^s_t\sL_r f(X_r)\dif r\ \mbox{ is a martingale}.\label{ERY1}
\end{align}
In other words, $\mP_{t,x}$ solves the martingale problem for $(\sL_t, C^2_b (\mR^d))$.

For any Borel set $E$, let
$$
  \sigma_E:=\inf\{s\geq 0:X_s\in E\},\qquad \tau_E:=\inf\{s\geq 0:X_s\notin E\},
$$
be the first hitting and exit time, respectively, of $E$.
Below for simplicity, we write
$$
 \cJ_\phi(t,x,y):=\frac{\kappa(t,x,y-x)}{ |y-x|^{d}\phi(|y-x|)}.
$$
We have the following L\'evy system of the Feller process $X$ (see \cite{Ch-Hu-Xi-Zh}).

\bl\label{lem:XLS}
  Let $f$ be a non-negative measurable function on $\mR_+\times\mR^d\times\mR^d$ that vanishes along the diagonal. 
  Then for every stopping time $T\geq t$,
  \begin{align}
    \mE_{t,x}\left[\sum_{t<r\leq T}f(r,X_{r-},X_r)\right]=\mE_{t,x}\left[\int^T_t\!\!\!\int_{\mR^d}
  f(r,X_r, y)\cJ_\phi(r,X_r,y)\dif y\dif r\right].     \label{sys}
  \end{align}
\el

We need the following two lemmas.
\bl\label{Le34}
There is a constant $\gamma_0\in(0,1)$ such that for all $\eps\in(0,1)$,
  \begin{align}
    \sup_{(t,x)\in\mR_+\times\mR^d}\mP_{t,x}\Big(\tau_{B(x,\eps)}\leq t+\gamma_0\phi(\eps)\Big)\leq {1}/{2}. \label{eqtaue}
  \end{align}
\el

\begin{proof}
  For simplicity, write $\tau:=\tau_{B(x,\eps)}$. By the strong Markov property of $X$, we have
  \begin{align}
    \mP_{t,x}\Big(\tau\leq t+r\Big)&\leq\mP_{t,x}\Big(\tau\leq t+r; X_{t+r}\in B(x,\eps/2)\Big) +\mP_{t,x}\Big(X_{t+r}\notin B(x,\eps/2)\Big)\no\\
    &=\mP_{t,x}\left(\mP_{\tau,X_\tau}\Big(X_{t+r}\in B(x,\eps/2)\Big); \tau\leq t+r\right)+\mP_{t,x}\Big(X_{t+r}\notin B(x,\eps/2)\Big)\no\\ &\leq\mP_{t,x}\left(\mP_{\tau,X_\tau}\Big(|X_{t+r}-X_\tau|\geq\eps/2\Big); \tau\leq t+r\right) +\mP_{t,x}\Big(X_{t+r}\notin B(x,\eps/2)\Big)\no\\
    &\leq 2\sup_{t\leq s\leq t+r}\sup_{x\in\mR^d}\mP_{s,x}\Big(|X_{t+r}-x|\geq\eps/2\Big),\label{AS}
  \end{align}
  where the second inequality is due to $|X_\tau-x|\geq \eps$ and $|X_{t+r}-x|\leq  {\eps}/2$.
  On the other hand, by \eqref{Tr} and the 
 heat kernel upper bound estimate in \eqref{GR1},  
  there is a constant $C>0$ such that for all $r\in(0,1)$, $t\leq s\leq t+r$ and $x\in\mR^d$,
  \begin{align*}
    \mP_{s,x}\Big(|X_{t+r}-x|\geq\eps/2\Big)
    &=\int_{{|x-y|\geq\eps/2}}p^\kappa_{s,t+r}(x, y)\dif y\leq   C(t+r-s)\int_{{|z|\geq  {\eps}/ 2}}
    \rho_\phi(t+r-s,z)\dif z\\ 
    & \leq Cr\int^\infty_{ {\eps}/2}\frac{\dif u}{u\phi(u)} 
     =Cr \int^\infty_{1}\frac{\dif u}{u\phi(\eps u/2)}
     \leq \frac{C r \cA^{(0)}_\phi  }{\phi(\eps/2)}\leq \frac{C_0r}{\phi(\eps)},
    \end{align*}
  where $\cA^{(0)}_\phi$ is defined
  in \eqref{Con0} with $i=0$. 
  Substituting this into \eqref{AS} yields
\begin{align}\label{DL4}
    \mP_{t,x}\Big(\tau_{B(x,\eps)}\leq t+r\Big)\leq \frac{C_0 r}{\phi(\eps)}.
\end{align}
Letting $r=\frac{\phi(\eps)}{2C_0}$ in \eqref{DL4}, we obtain \eqref{eqtaue}  with $\gamma_0=\frac1{2C_0}$.
\end{proof}

\bl\label{Le35}
 Let $\gamma_0$ be  the constant from Lemma \ref{Le34}. 
  For all $\gamma\in(0,\gamma_0]$,
 there exists a constant $c_1>0$ such that for all $t>0$, $\eps\in(0,1)$ and $x,y\in\mR^d$ with $|x-y|\geq 2\eps$,
  \begin{align}\label{GF1}
    \mP_{t,x}\Big(\sigma_{B(y,\eps)}<t+\gamma \phi(\eps)\Big)\geq c_1\frac{\eps^{d} \phi(\eps)}{|x-y|^d\phi(|x-y|)}.
  \end{align}
\el

\begin{proof}
  For $\eps\in(0,1)$ and $\gamma\in(0,\gamma_0]$, by (\ref{eqtaue}) we have
  \begin{align}\label{CC2}
  \mE_{t,x} \left[  \left(t+\gamma \phi(\eps) \right)\wedge \tau_{B(x,\eps) } -t \right] 
  \geq \gamma \phi(\eps)\mP_{t,x}
    \Big(\tau_{B(x,\eps)}\geq t+\gamma \phi(\eps)\Big) \geq \frac{\gamma \phi(\eps)}{2}.
  \end{align}
  Noticing that under $\mP_{t,x}$, 
  $$
    X_r\notin B(y,\eps)\quad \text{when } \  t<r<(t+\gamma \phi(\eps))\wedge \tau_{B(x,\eps)},
  $$
  we have
  $$
     {\1}_{X_{(t+\gamma \phi(\eps))\wedge \tau_{B(x,\eps)}}\in B(y,\eps)}=\sum_{t<r\leq (t+\gamma \phi(\eps))
    \wedge \tau_{B(x,\eps)}} {\1}_{X_{r}\in B(y,\eps)}.
  $$
 By the L\'evy system formula 
 (\ref{sys}) and the definition of $\cJ_\phi$, we have
  \begin{align}
    \mP_{t,x}\Big(\sigma_{B(y,\eps)}<t+\gamma \phi(\eps)\Big)&\geq\mP_{t,x}\Big(X_{(t+\gamma \phi(\eps))\wedge \tau_{B(x,\eps)}}\in B(y,\eps)\Big)\no\\
    &=\mE_{t,x}\int_t^{(t+\gamma \phi(\eps))\wedge \tau_{B(x,\eps)}}\!\!\int_{B(y,\eps)}
    \frac{\kappa(t,X_r,z-X_r)}{|z-X_r|^d\phi(|z-X_r|)}\dif z\dif r.\label{CC3}
  \end{align}
  Since $|x-y|\geq 2\eps$, we have for all $z\in B(y,\eps)$ and $X_r\in B(x,\eps)$,
  $$
    |z-X_r|\leq|y-z|+|x-y|+|X_r-x| < 2|x-y|.
  $$
  Thus  by \eqref{CC3} and \eqref{CC2}, we have
  \begin{align*}
    \mP_{t,x}\Big(\sigma_{B(y,\eps)}<t+\gamma \phi(\eps)\Big)&\geq \frac{\gamma \phi(\eps)}{2}\int_{B(y,\eps)}\frac{\kappa_0^{-1}  }
    {(2|x-y|)^d\phi(2|x-y|)} \dif z 
    \geq c_2\frac{\eps^{d} \phi(\eps)}{|x-y|^d\phi(|x-y|)}.
  \end{align*}
  This proves the lemma.
\end{proof}

Now we can give

\begin{proof}[Proof of lower bound in Theorem \ref{heat}{\bf (a)}]
 Let $\delta >0$ be the constant 
  in Corollary \ref{Co42}. 
 We claim that by \eqref{DL8},
 for any   $0\leq t<s\leq T$, $x, y\in \R^d$  and  
 $n\in\mN$, there is a constant $C_n>0$ such that
  \begin{align}\label{ER2}
    p^\kappa_{t,s}(x,y)\geq C_n\phi^{-1}(s-t)^{-d} \quad \hbox{whenever  } |x-y|\leq\phi^{-1}((s-t)/2^n)\leq \delta.
  \end{align}
  Indeed, if 
  $|x-y|\leq \phi^{-1}((s-t)/2)\leq \delta$, 
  then by the Chapman-Kolmogorov equation, 
  {\setlength{\arraycolsep}{0.2pt}
  \begin{eqnarray*}
    p^\kappa_{t,s}(x,y)&=&\int_{\mR^d}p^\kappa_{t,\frac{t+s}{2}}(x,z)p^\kappa_{\frac{t+s}{2},s}(z,y)\dif z\\
    &\geq&\int_{B(\frac{x+y}{2}, \phi^{-1}((s-t)/2 ))}
    p^\kappa_{t,\frac{t+s}{2}}(x,z)p^\kappa_{\frac{t+s}{2},s}(z,y)\dif z\\
    &\stackrel{(\ref{DL8})}{\geq}& c^2_1\phi^{-1}((s-t)/2)^{-2d} \, \Vol\left(B \left(\tfrac{x+y}{2},
    \phi^{-1}((s-t)/2) \right)\right)\\  
    &\stackrel{(\ref{e:1.1}) }{\succeq} & \phi^{-1}(s-t)^{-d}.   
  \end{eqnarray*}}
Iterating the above estimates establishes  the claim \eqref{ER2}. 

  Now fix $T>0$ and choose $n$ large enough so that 
$$
T/2^n\leq \phi (\delta) ,  \quad\hbox{or equivalent, }  \  \phi^{-1}(T/2^n) \leq \delta.
$$ 
Consider $x, y\in \R^d$ and $0\leq t<s\leq T$.  
  If $|x-y|\leq\phi^{-1}((s-t)/2^n)$, 
  we immediately get from  \eqref{ER2} the lower bound  for $p^\kappa_{t, s}(x, y)$  in  \eqref{GR1}.
  It remains to consider the case that $|x-y|>\phi^{-1}((s-t)/2^n)$. 
  Define 
 \begin{equation}\label{e:4.56}
  \eps = \tfrac13 \phi^{-1} \left((s-t)/2^{n+1} \right) \quad \hbox{so } \  (s-t)/2^{n } = 2 \phi (3\eps).
  \end{equation} 
     Let $\gamma_0\in(0,1)$ be the constant in Lemma \ref{Le34}. 
  By the strong Markov property of $X$ and Lemma \ref{Le35}, we have for any 
  $|x-y|\geq 3\eps$,
  {\setlength{\arraycolsep}{0pt}
  \begin{eqnarray*}
    \mP_{t,x}\Big(X_{t+2\gamma_0 \phi(\eps)}\in B\big(y,2\eps\big)\Big)&\geq& \mP_{t,x}\left(\sigma:=\sigma_{B(y,\eps)}<t+\gamma_0 \phi(\eps);
    \sup_{s\in[\sigma,\sigma+\gamma_0 \phi(\eps)]}|X_s-X_\sigma|< \eps\right)\no\\
    &=&\mE_{t,x}\left[\mP_{\sigma,X_\sigma}\left(\sup_{s\in[\sigma,\sigma+\gamma_0 \phi(\eps)]}|X_s-X_\sigma| < \eps\right);
    \sigma_{B(y,\eps)} < t+\gamma_0 \phi(\eps)\right]\no\\
    &\geq& \inf_{r,z}\mP_{r,z}\left(\tau_{B(z,\eps)}> r+\gamma_0 \phi(\eps)\right) \mP_{t,x}\left(\sigma_{B(y,\eps)}< t+\gamma_0 \phi(\eps)\right)\no\\
    &\stackrel{(\ref{eqtaue})}{\geq}&\frac{1}{2} \mP_{t,x}\left(\sigma_{B(y,\eps)}<t+\gamma_0 \phi(\eps)\right) \\
    & \stackrel{(\ref{GF1})}{\geq} & \frac{\eps^{d} \phi(\eps)}{|x-y|^d\phi(|x-y|)}.
  \end{eqnarray*}}
  Hence we have for any $x, y\in \R^d$ with 
  $|x-y|\geq 3\eps$,
  \begin{eqnarray*}
    p^\kappa_{t,s}(x,y
    &\geq & \int_{B(y,2\eps)}p^\kappa_{t,t+2\gamma_0\phi(\eps)}\big(x,z\big) 
    p^\kappa_{t+2\gamma_0 \phi(\eps),s}(z,y)\dif z\\
    &\geq & \inf_{z\in B(y,2\eps)}p^\kappa_{t+2\gamma_0 \phi(\eps),s}(z,y)\mP_{t,x}\Big(X_{t+2\gamma_0 \phi(\eps)}\in B\big(y,2\eps\big)\Big)\\
  & \stackrel{(\ref{ER2})}{\geq} & C_1
    \phi^{-1}(s-t)^{-d}\cdot\frac{\eps^{d} \phi(\eps)}{|x-y|^d\phi(|x-y|)} \\
 &    \stackrel{(\ref{e:1.1})}{\geq}  & C_2  (s-t)\rho_\phi(s-t,x-y),
  \end{eqnarray*}
where in the third inequality when apply \eqref{ER2} we used the fact that from \eqref{e:4.56}
$$ 
\phi (\delta ) > \frac{s-(t+2\gamma_0 \phi(\eps))}{2^n}\geq 2\phi (3\eps) -\gamma_0 \phi(\eps)\geq \phi (2\eps) > \phi (|z-y|)
\quad \hbox{for } z\in B(y,2\eps). 
$$
This establishes  the lower bound  for $p^\kappa_{t, s}(x, y)$ in  \eqref{GR1}  on $\mD^T_0$. 
\end{proof}

\section{Examples}\label{S:5}

In this section, 
we discuss the assumptions \eqref{e:1.1}, 
 \eqref{Con0}
and \eqref{Con3}, and give some examples 
 that satisfy these conditions and therefore our main results apply.

\bx\rm
Let $\phi(r)=r^\alpha$ with $\alpha\in(0,2)$. It is easy to see that \eqref{e:1.1}, 
 \eqref{Con0}
and \eqref{Con3} hold.
In particular, Theorem \ref{heat} in this case extends the main results in \cite{Ch-Zh} and  \cite{Ch-Zh2} to time-dependent and Dini's continuous kernels $\kappa(t,x,z)$.
\ex

\bx\rm
Let $0<\beta_1\leq\beta_2<2$ and $\phi$ an increasing function on $[0, \infty)$ so that there are positive constants
$c^\phi_2\geq c^\phi_1$ such that  
$$
c^\phi_1\left(\frac{R}{r}\right)^{\beta_1} \leq\frac{\phi  (R)}{\phi  (r)}\leq c^\phi_2\left(\frac{R}{r}\right)^{\beta_2}
\quad \hbox{for any }  0<r<R<\infty .
$$
This is the case, for example, when 
$$ \phi (r) = \int_{\beta_1}^{\beta_2} r^\alpha \nu (\dif \alpha) \quad \hbox{or} \quad 
\phi (r) = 1 / {\int_{\beta_1}^{\beta_2} r^{-\alpha} \nu (\dif \alpha)},
$$
where $\nu$ is a probability measure on $[\beta_1, \beta_2]$.   
Clearly, $\phi$ satisfies \eqref{e:1.1} and 
 ({\bf A$^{(0)}_\phi$}),
as well as 
 ({\bf A$^{(1)}_\phi$})
of  ${\rm Case^\phi_1}$ and ${\rm Case^\phi_2}$.
Property  
 ({\bf A$^{(1)}_\phi$})
holds when $\beta_1>1$, under which   ${\rm Case^\phi_3}$ occurs.  
When $\kappa(t,x,z)$ is symmetric in $z$,  
Theorem \ref{heat} in particular extends the main results in 
\cite{Ch-Zh, Ch-Zh2, KSV}. See  (iii) and (iv) of Remark \ref{R:1.4}  for the gradient estimate.

\ex
\bx \label{E:5.3} \rm
 Let $\phi(r)=r^{\alpha}\1_{\{r\leq 1\}}+r^\beta \1_{\{r>1\}}$ 
 with $\alpha>0$ and $\beta >0$. 
 Note that \eqref{e:1.1} holds with $\beta_1:=\alpha$ and $\beta_2:=\alpha \vee \beta$.
One can check that ({\bf A$^{(0)}_\phi$}) is satisfied for any $\alpha\in(0,2)$ and  $\beta>0$,
while ({\bf A$^{(1)}_\phi$})
holds for all   $\alpha\in(0, 2)$ and  $\beta>0$ except 
 when $\alpha =1$ and $\beta >1$ (corresponding to Case$^\phi_3$).
As we noted in \eqref{e:1.13}, for $\ell (r)=r^\eta$ on $(0, 1]$,  condition 
\eqref{Con3} holds if and only if $\alpha>1/2$ and $\alpha+\eta>1$. 
Note that the function $\phi$ can have  any polynomial growth as $r\to\infty$.
 When  $\ell \in \sD_0$,   condition \eqref{Con3} holds if $\alpha \geq 1$ 
with 
$$
 M^\phi_\ell (t) =\frac{\alpha}{\ell (t)} \int_0^t \frac{\ell (r)}r  r^{\alpha-1} \dif r + \frac{\alpha}{2\alpha-1} t^{\alpha-1}
  \quad \hbox{on } (0, 1]
$$
When $\alpha >1$,  clearly 
$ M^\phi_\ell (t) \asymp t^{\alpha-1} $ on $(0, 1]$ and thus  the gradient estimate \eqref{Grad} takes the form
 \begin{equation}\label{e:5.1A}
 |\nabla p^\kappa_{t,s}(\cdot,y)(x)|\leq 
   \wt   c_3  (s-t)^{1-(1/\alpha)}  
\rho_\phi(s-t,x-y)\quad  \hbox{on } \mD^T_0.
 \end{equation} 
 When $\alpha =1$,  by the same calculation as that for \eqref{e:1.31}, 
we have  
$ M^\phi_\ell (t) \asymp  \frac{\Gamma_\ell (t)}{\ell (t)} $ on $(0, 1]$ 
and so  the gradient estimate \eqref{Grad} has the form
  \begin{align}\label{e:5.1B}
 |\nabla p^\kappa_{t,s}(\cdot,y)(x)|\leq 
   \wt   c_3   \frac{\Gamma_\ell (s-t)}{\ell (s-t)}  
\rho_\phi(s-t,x-y)\quad  \hbox{on } \mD^T_0.
\end{align}
 \ex

\smallskip

We have the following   more general result.
\bp
  Suppose that $\phi$ is an increasing function on $\mR_+$ so that $\phi\in\sR_\alpha$ on $(0, 1]$ 
 and $r\mapsto\phi(1/r)\in\sR_\beta$  for $\alpha \in (0, 1)$ and
 $\beta <0$
 or  for $ \alpha \in (1, 2)$ and $\beta \in (-\infty, -1) \cup (-1, 0)$. 
Then \eqref{e:1.1} and 
 ({\bf A$^{(1)}_\phi$}) hold
\ep
\begin{proof}
By definition, there are some slowly varying functions $\ell,\ell'\in\sS_0$ so that 
$$
\phi(r)=r^\alpha\ell(r) \hbox{ when  } r\in(0,1)  \quad \hbox{and} \quad 
\phi(r)=r^{- \beta}\ell'(1/r) 
\hbox{ when }  r\in(1,\infty).
$$
Take $0<\eps\ < \alpha \wedge (-\beta)$.  By \eqref{e:3.1}, we see that \eqref{e:1.1} 
holds with $\beta_1:= \alpha -\eps$ and $\beta_2:= (\alpha +\eps)\vee (\eps-\beta)$. 

Suppose $\alpha\in(0,1)$ and $\beta <0$. Then Case$^\phi_1$ holds. 
By \eqref{e:3.1} again, we have 
 for $\delta\in(0,1-\alpha)$,
$$
\sup_{\lambda \in (0, 1]} \int^1_0\frac{\gamma_\phi (r) \phi(\lambda)}{r \phi(r\lambda)}\dif r =
\sup_{\lambda \in (0, 1]} \int^1_0\frac{\phi(\lambda)}{\phi(r\lambda)}\dif r\lesssim 
 \int^1_0r^{-\delta-\alpha}\dif r<\infty,
$$
and for $\delta\in(0, -\beta)$
\begin{eqnarray}
 \sup_{\lambda \in (0, 1]} \int_1^\infty \frac{\gamma_\phi (r) \phi(\lambda)}{r \phi(r\lambda)}\dif r
& =&   \sup_{\lambda \in (0, 1]} \int^\infty_1\frac{\phi(\lambda)}{r\phi(r\lambda)}\dif r 
  \nonumber \\
&=& \sup_{\lambda \in (0, 1]} \phi (\lambda)   \left( \int_{\lambda}^1 \frac{1}{s\phi (s)} \dif s
+ \int_1^\infty  \frac{1}{s\phi (s)} \dif s \right)  \nonumber \\
&  \lesssim &  \sup_{\lambda \in (0, 1]} \left( 1+ \phi (\lambda) \int^\infty_1
\frac{\dif r}{s^{1-\beta-\delta}} \dif s \right) <\infty,    \label{e:5.1a}
\end{eqnarray}
where in the last inequality we used \eqref{e:1.6}. 
Thus  
({\bf A$^{(1)}_\phi$})
holds when $\alpha\in(0,1)$ and $\beta <0$.

When $\alpha\in(1,2)$, we have by 
 \eqref{e:3.1} 
that $\int_{0+} \frac1{\phi (r)} \dif r =\infty$
and for 
 $\delta\in(0,2-\alpha)$, 
$$
\sup_{\lambda \in (0, 1]} \int^1_0\frac{\gamma_\phi (r) \phi(\lambda)}{r \phi(r\lambda)}\dif r
= \sup_{\lambda \in (0, 1]} \int^1_0\frac{r^2 \phi(\lambda)}{r \phi(r\lambda)}\dif r
\lesssim 
 \int_0^1 r^{1-\delta -\alpha} \dif r <\infty.
$$
When $\beta \in (-1, 0)$, Case$^\phi_2$ holds. In this case,  
 we have by  \eqref{e:5.1a}, 
$$
 \sup_{\lambda \in (0, 1]} \int_1^\infty \frac{\gamma_\phi (r) \phi(\lambda)}{r \phi(r\lambda)}\dif r =
  \sup_{\lambda \in (0, 1]} \int_1^\infty\frac{   \phi(\lambda)}{r \phi(r\lambda)}\dif r
  <\infty. 
$$
When $\beta \in (-\infty, -1)$, Case$^\phi_3$ holds. In this case,   
  for some $0< \delta <   \min\{ \alpha -1,  -\beta -1\}$,  we have by \eqref{e:3.1} that 
\begin{eqnarray*}
 \sup_{\lambda \in (0, 1]} \int_1^\infty \frac{\gamma_\phi (r) \phi(\lambda)}{r \phi(r\lambda)}\dif r 
 &=& \sup_{\lambda \in (0, 1]} \int_1^\infty\frac{   \phi(\lambda)}{ \phi(r\lambda)}\dif r
 =  \sup_{\lambda \in (0, 1]}  \frac{1}{\lambda} \int_{\lambda}^\infty \frac{\phi (\lambda )}{\phi (s)} \dif s \\
&=&  \sup_{\lambda \in (0, 1]}     \left( \frac{1}{\lambda}  \int_{\lambda}^1 \frac{  \phi (\lambda )}{\phi (s)} \dif s
+    \frac{\phi (\lambda )}{\lambda}    \int_1^\infty  \frac{1}{\phi (s)} \dif s \right) \\
&\lesssim & \sup_{\lambda \in (0, 1]}    \left( \frac{1}{\lambda}  \int_{\lambda}^1 (\lambda /s)^{\alpha -\delta} \dif s
+   \lambda^{\alpha -\delta -1}   \int_1^\infty  s^{\delta  +\beta}  \dif s \right) \\
&\lesssim & 1+1<\infty. 
\end{eqnarray*}
Hence  ({\bf A$^{(1)}_\phi$})
holds when $ \alpha \in (1, 2)$ and  $\beta \in (-\infty, -1) \cup (-1, 0)$. 
\end{proof}

\bx\rm
Let $\ell(s)=(\log (1/s) )^{-2}$ and $\phi(s)=s\log (1/s) $ for $ 0< s \ll  1$. It is easy to see that
$$
\ell\in\sS_0\cap\sD_0 \quad \hbox{and} \quad \int_{0+}1/\phi(s)\dif s=\infty.
$$
However,
$$
 \int_{0+}\frac{\ell(s)}{s}\dif \phi(s)=\int_{0+}\frac{(\log (1/s) +1)}{s(\log (1/s) )^2}\dif s=\infty.
$$
Hence in this case,  \eqref{Con3} does not hold.
\ex

\medskip

\bx \label{E:5.6} \rm
Suppose that $\ell \in \sR_\eta$ for some $\eta\geq 0$ and $\phi$ is an   increasing function on $[0, 1]$ such that there
are $\beta_1>0$ and $c_1>0$ so that 
\begin{equation}\label{e:5.1} 
\phi (R)/\phi (r) \geq c_1 (R/r)^{\beta_1} \quad \hbox{for  any } 0<r<R\leq 1.
\end{equation}  
Then condition \eqref{Con3} holds if  $\beta_1>1/2$ and $\beta_1 + \eta >1$,  and in this case
\begin{equation} \label{e:5.2} 
M^\phi_\ell (t)   
  :=\int^t_0\frac{1}{r}\left(\frac{\ell(r)}{\ell(t)}+\frac{\phi(r)}{\phi(t)}\right)\dif \phi(r) \asymp
\frac{\phi (t)}t      \quad \hbox{ for } t\in (0, 1].
\end{equation} 
Indeed,  there is $\ell_0 \in \sS_0$ so that $\ell (r)=r^\eta \ell_0(r)$.
By Proposition \ref{Pr32}(i), 
for any $\delta\in(0,(\beta_1 + \eta -1) \wedge 1)$, there is  $c_0>0$ so that 
$$
\frac{\ell (s)}{\ell (t)} \leq c_0 \left(\frac{s}{t}\right)^{\eta -\delta}  \quad \hbox{for } 0<s<t\leq 1.
$$
Thus  by integration by parts and the fact that  $\beta_1+\eta -\delta >1$,  we have for $t\in (0, 1]$, 
\begin{equation}\label{e:5.3}
\int_0^t \frac{\ell (r)}{r \ell (t)} \dif \phi (r)
\lesssim   t^{\delta -\eta } \int_0^t r^{\eta -\delta -1} \dif \phi (r) 
=  t^{\delta -\eta} \left( t^{\eta-\delta - 1} \phi (t) + (1-\eta-\delta ) \int_0^t \phi (r) r^{\eta-\delta -2} \dif r\right) . 
\end{equation} 
By \eqref{e:5.1},
$$
\int_0^t \phi (r) r^{\eta -\delta -2} \dif r
= \phi (t) \int_0^t \frac{\phi (r)}{\phi (t)}  r^{\eta -\delta -2} \dif r
 \leq  \frac{\phi (t)} {c_1} \int_0^t (r/t)^{\beta_1} r^{\eta -\delta -2} \dif r 
 = \frac{ \phi (t) t^{\eta -\delta - 1}}{c_1 ( \beta_1 + \eta -\delta  -1 ) }. 
$$
This together with \eqref{e:5.3} shows that 
\begin{equation} \label{e:5.4}
\int_0^t \frac{\ell (r)}{r \ell (t)} \dif \phi (r)  \asymp  \frac{\phi (t)}t   \quad \hbox{ for } t\in (0, 1].
\end{equation} 
On the other hand, for $t\in (0, 1]$, by integration by parts and the assumption that $\beta_1 > 1/2$, 
\begin{equation}\label{e:5.5} 
 \int_0^t \frac{\phi (r)}{r \phi (t)} \dif \phi (r)
  =  \frac1{2\phi (t)} \int_0^t r^{-1} \dif (\phi (r)^2) = \frac1{2\phi (t)} 
 \left( t^{-1} \phi (t)^2 + \int_0^t r^{-2} \phi (r)^2   \dif r \right) .
 \end{equation}
 By \eqref{e:5.1}, 
 $$
 \int_0^t r^{-2} \phi (r)^2   \dif r =     \phi (t)^2 \int_0^t r^{-2} \left(\frac{\phi (r)}{\phi (t)} \right)^2   \dif r
 \leq \frac{\phi (t)^2}{c_1} \int_0^t r^{-2} (r/t)^{2\beta_1} \dif r 
 = \frac{\phi (t)^2}{c_1 (2\beta_1-1) t}. 
$$
It follows from \eqref{e:5.5} that 
$$
 \int_0^t \frac{\phi (r)}{r \phi (t)} \dif \phi (r) 
 \asymp 
 \frac{\phi (t)}t   \quad \hbox{ for } t\in (0, 1].
$$
This together with \eqref{e:5.4} proves the claim \eqref{e:5.2}. 
\ex

\bigskip

\end{document}